\documentclass{gtart}

\def\ifplaintex{\expandafter\ifx\csname documentclass\endcsname\relax}

\ifplaintex 
\else
\fi

\expandafter\ifx\csname beginpicture\endcsname\relax
\expandafter\ifx\csname documentclass\endcsname\relax
\input pictex \else
\input prepictex \input pictex \input postpictex \fi\fi

\def\gt{{\mathsurround=0pt\it $\cal G\mskip-2mu$eometry \&\ 
$\cal T\!\!$opology}}        

\def\gtp{{\mathsurround=0pt\it $\cal G\mskip-2mu$eometry \&\ 
$\cal T\!\!$opology $\cal P\!$ublications}}  

\def\lognumber#1{\def\thelognumber{#1}}
\def\volumenumber#1{\def\thevolumenumber{#1}}
\def\papernumber#1{\def\thepapernumber{#1}}
\def\volumeyear#1{\def\thevolumeyear{#1}}

\def\pagenumbers#1#2{\def\startpage{#1}\def\finishpage{#2}}
\def\published#1{\def\publishdate{#1}}
\def\proposed#1{\def\theproposer{#1}}
\def\seconded#1{\def\theseconders{#1}}
\def\received#1{\def\receiveddate{#1}}
\def\revised#1{\def\reviseddate{#1}}
\def\accepted#1{\def\accepteddate{#1}}

\let\\\par\let\thelognumber\relax
\let\thevolumenumber\relax\let\thepapernumber\relax
\let\thevolumeyear\relax\let\thesamplenumber\relax\let\startpage\relax
\let\finishpage\relax\let\publishdate\relax\let\receiveddate\relax
\let\reviseddate\relax\let\accepteddate\relax\let\theasciititle\relax
\let\theasciiauthors\relax
\let\theasciiabstract\relax
\let\theasciiemail\relax\let\theshortauthors\relax\let\theshorttitle\relax

\long\def\maketitlep{   

\count0=\startpage

\gt\hfill      
\beginpicture
\setcoordinatesystem units <0.33truein, 0.33truein> point at 2.2 0.9
\setplotsymbol ({$\cal G$})
\plotsymbolspacing=9truept
\circulararc 315 degrees from 0 1 center at 0 0
\setplotsymbol ({$\cal T$})
\circulararc 315 degrees from 1 -1 center at 1 0
\endpicture
\break
{\small\ifx\thesamplenumber\relax 
Volume \else Sample
\fi\thevolumenumber\ (\thevolumeyear)
\startpage--\finishpage\nl
Published: \publishdate}
\vglue 0.5truein plus 0.4fil minus 0.1truein

{\parskip=0pt\leftskip 0pt plus 1fil\def\\{\par\smallskip}{\ifplaintex\large
\else\Large\fi\bf\thetitle}\par\medskip}   

\vglue 0pt plus 0.1fil 

{\parskip=0pt\leftskip 0pt plus 1fil\def\\{\par}{\sc\theauthors}
\par\medskip}

\vglue 0pt plus 0.1fil 

{\small\parskip=0pt\let\newline\\
{\leftskip 0pt plus 1fil\def\\{\par}{\sl\theaddress}\par}
\expandafter\ifx\theemail\relax    
\relax\else\vglue 5pt plus 0.02fil minus 2pt\def\\{\stdspace{\rm 
and}\stdspace} 
\cl{Email:\stdspace\tt\theemail}\fi
\ifx\theurl\relax                  
\relax\else\vglue 5pt plus 0.02fil minus 2pt\def\\{\stdspace{\rm 
and}\stdspace}
\cl{URL:\stdspace\tt\theurl}\fi\par}

\vglue 7pt plus 0.3fil minus 3pt

{\bf Abstract}
\vglue 5pt plus 0.1fil minus 2pt

\theabstract

\vglue 7pt plus 0.3fil minus 3pt

{\bf AMS Classification numbers}\quad Primary:\quad \theprimaryclass

Secondary:\quad \thesecondaryclass

\vglue 5pt plus 0.3fil minus 2pt

{\bf Keywords}\quad \thekeywords

\vglue 10pt plus 0.5fil minus 5pt

{\small  Proposed: \theproposer\hfill Received: \receiveddate\nl
Seconded: \theseconders\hfill 
\ifx\reviseddate\relax                         
Accepted: \accepteddate                        
\else
Revised: \reviseddate                          
\fi}
\eject
}       

\let\maketitlepage\maketitlep
\let\maketitle\maketitlepage

\font\phead=cmsl9 scaled 950
\font=cmsl9 scaled 1050
\font\pnum=cmbx10 scaled 913
\font\lnum=cmbx10 
\font\pfoot=cmsl9 scaled 950
\font=cmsl9 scaled 1050
\ifplaintex
\headline{\vbox to 0pt{\vskip -4.5mm\line{\small\phead\ifnum
\count0=\startpage ISSN 1364-0380 (on line)
1465-3060 (printed) \hfill {\pnum\folio}\else\ifodd\count0\def\\{ }%
\ifx\theshorttitle\relax\thetitle\else\theshorttitle\fi\hfill{\pnum\folio}
\else\def\\{ and }{\pnum\folio}\hfill\ifx\theshortauthors\relax\theauthors
\else\theshortauthors\fi\fi\fi}\vss}}
\footline{\vbox to 0pt{\vglue 0mm\line{\small\pfoot\ifnum\count0=\startpage
\copyright\ \gtp\hfill\else
\gt, Volume \thevolumenumber\ (\thevolumeyear)\hfill\fi}\vss
}}
\else
\makeatletter
\def\@oddhead{{\small\lhead\ifnum\count0=\startpage ISSN 1364-0380 (on line)
1465-3060 (printed) \hfill {\lnum\number\count0}\else\ifodd\count0
\def\\{ }\ifx\theshorttitle\relax \thetitle \else\theshorttitle\fi\hfill
{\lnum\number\count0}\else\def\\{ and }{\lnum\number\count0}
\hfill\ifx\theshortauthors\relax 
\theauthors\else\theshortauthors\fi\fi\fi}}\def\@evenhead{\@oddhead}
\def\@oddfoot{\small\lfoot\ifnum\count0=\startpage\copyright\ \gtp\hfill\else
\gt, Volume \thevolumenumber\ (\thevolumeyear)\hfill\fi}
\def\@evenfoot{\@oddfoot}
\makeatother
\fi

\newwrite\gtoutfile
\long\gdef\makeheadfile{  
{\def\\{, }\def\s{ }
\immediate\openout\gtoutfile head.xxx
\immediate\write\gtoutfile{To: math@arxiv.org}
\immediate\write\gtoutfile{Subject: put or rep NNNNN:pppp}
\immediate\write\gtoutfile{--text follows this line--}
\immediate\write\gtoutfile{Proxy-for: \ifx\theasciiauthors\relax
\theauthors\else\theasciiauthors\fi\s<\ifx\theasciiemail\relax\theemail\else\theasciiemail\fi>}
\immediate\write\gtoutfile{\noexpand\\}
\immediate\write\gtoutfile{Authors: \ifx\theasciiauthors\relax
\theauthors\else\theasciiauthors\fi}
\immediate\write\gtoutfile{Title: \ifx\theasciititle\relax
\thetitle\else\theasciititle\fi}
\immediate\write\gtoutfile{Subj-class: GT or SG or MG etc}
\immediate\write\gtoutfile{MSC-class: \theprimaryclass\ifx\thesecondaryclass\relax\else, \thesecondaryclass\fi}
\immediate\write\gtoutfile{Journal-ref: Geom. Topol. \thevolumenumber
(\thevolumeyear) \startpage-\finishpage}
\immediate\write\gtoutfile{Comments: Published by Geometry and Topology at}
\immediate\write\gtoutfile{\s\s http://www.maths.warwick.ac.uk/gt/GTVol\thevolumenumber/paper\thepapernumber.abs.html}
\immediate\write\gtoutfile{\noexpand\\}
\immediate\write\gtoutfile{}
\ifx\theasciiabstract\relax
\immediate\write\gtoutfile{\theabstract}\else
\immediate\write\gtoutfile{\theasciiabstract}\fi
\immediate\write\gtoutfile{}
\immediate\write\gtoutfile{\noexpand\\}
\immediate\write\gtoutfile{}
\immediate\closeout\gtoutfile}}  

\def\maketitlepage{\maketitlep\makeheadfile}
\let\maketitle\maketitlepage

\def\ifplaintex{\expandafter\ifx\csname documentclass\endcsname\relax}

\ifplaintex 
\else
\fi

\expandafter\ifx\csname beginpicture\endcsname\relax
\expandafter\ifx\csname documentclass\endcsname\relax
\input pictex \else
\input prepictex \input pictex \input postpictex \fi\fi

\def\gt{{\mathsurround=0pt\it $\cal G\mskip-2mu$eometry \&\ 
$\cal T\!\!$opology}}        

\def\gtp{{\mathsurround=0pt\it $\cal G\mskip-2mu$eometry \&\ 
$\cal T\!\!$opology $\cal P\!$ublications}}  

\def\lognumber#1{\def\thelognumber{#1}}
\def\volumenumber#1{\def\thevolumenumber{#1}}
\def\papernumber#1{\def\thepapernumber{#1}}
\def\volumeyear#1{\def\thevolumeyear{#1}}

\def\pagenumbers#1#2{\def\startpage{#1}\def\finishpage{#2}}
\def\published#1{\def\publishdate{#1}}
\def\proposed#1{\def\theproposer{#1}}
\def\seconded#1{\def\theseconders{#1}}
\def\received#1{\def\receiveddate{#1}}
\def\revised#1{\def\reviseddate{#1}}
\def\accepted#1{\def\accepteddate{#1}}

\let\\\par\let\thelognumber\relax
\let\thevolumenumber\relax\let\thepapernumber\relax
\let\thevolumeyear\relax\let\thesamplenumber\relax\let\startpage\relax
\let\finishpage\relax\let\publishdate\relax\let\receiveddate\relax
\let\reviseddate\relax\let\accepteddate\relax\let\theasciititle\relax
\let\theasciiauthors\relax
\let\theasciiabstract\relax
\let\theasciiemail\relax\let\theshortauthors\relax\let\theshorttitle\relax

\long\def\maketitlep{   

\count0=\startpage

\gt\hfill      
\beginpicture
\setcoordinatesystem units <0.33truein, 0.33truein> point at 2.2 0.9
\setplotsymbol ({$\cal G$})
\plotsymbolspacing=9truept
\circulararc 315 degrees from 0 1 center at 0 0
\setplotsymbol ({$\cal T$})
\circulararc 315 degrees from 1 -1 center at 1 0
\endpicture
\break
{\small\ifx\thesamplenumber\relax 
Volume \else Sample
\fi\thevolumenumber\ (\thevolumeyear)
\startpage--\finishpage\nl
Published: \publishdate}
\vglue 0.5truein plus 0.4fil minus 0.1truein

{\parskip=0pt\leftskip 0pt plus 1fil\def\\{\par\smallskip}{\ifplaintex\large
\else\Large\fi\bf\thetitle}\par\medskip}   

\vglue 0pt plus 0.1fil 

{\parskip=0pt\leftskip 0pt plus 1fil\def\\{\par}{\sc\theauthors}
\par\medskip}

\vglue 0pt plus 0.1fil 

{\small\parskip=0pt\let\newline\\
{\leftskip 0pt plus 1fil\def\\{\par}{\sl\theaddress}\par}
\expandafter\ifx\theemail\relax    
\relax\else\vglue 5pt plus 0.02fil minus 2pt\def\\{\stdspace{\rm 
and}\stdspace} 
\cl{Email:\stdspace\tt\theemail}\fi
\ifx\theurl\relax                  
\relax\else\vglue 5pt plus 0.02fil minus 2pt\def\\{\stdspace{\rm 
and}\stdspace}
\cl{URL:\stdspace\tt\theurl}\fi\par}

\vglue 7pt plus 0.3fil minus 3pt

{\bf Abstract}
\vglue 5pt plus 0.1fil minus 2pt

\theabstract

\vglue 7pt plus 0.3fil minus 3pt

{\bf AMS Classification numbers}\quad Primary:\quad \theprimaryclass

Secondary:\quad \thesecondaryclass

\vglue 5pt plus 0.3fil minus 2pt

{\bf Keywords}\quad \thekeywords

\vglue 10pt plus 0.5fil minus 5pt

{\small  Proposed: \theproposer\hfill Received: \receiveddate\nl
Seconded: \theseconders\hfill 
\ifx\reviseddate\relax                         
Accepted: \accepteddate                        
\else
Revised: \reviseddate                          
\fi}
\eject
}       

\let\maketitlepage\maketitlep
\let\maketitle\maketitlepage

\font\phead=cmsl9 scaled 950
\font=cmsl9 scaled 1050
\font\pnum=cmbx10 scaled 913
\font\lnum=cmbx10 
\font\pfoot=cmsl9 scaled 950
\font=cmsl9 scaled 1050
\ifplaintex
\headline{\vbox to 0pt{\vskip -4.5mm\line{\small\phead\ifnum
\count0=\startpage ISSN 1364-0380 (on line)
1465-3060 (printed) \hfill {\pnum\folio}\else\ifodd\count0\def\\{ }%
\ifx\theshorttitle\relax\thetitle\else\theshorttitle\fi\hfill{\pnum\folio}
\else\def\\{ and }{\pnum\folio}\hfill\ifx\theshortauthors\relax\theauthors
\else\theshortauthors\fi\fi\fi}\vss}}
\footline{\vbox to 0pt{\vglue 0mm\line{\small\pfoot\ifnum\count0=\startpage
\copyright\ \gtp\hfill\else
\gt, Volume \thevolumenumber\ (\thevolumeyear)\hfill\fi}\vss
}}
\else
\makeatletter
\def\@oddhead{{\small\lhead\ifnum\count0=\startpage ISSN 1364-0380 (on line)
1465-3060 (printed) \hfill {\lnum\number\count0}\else\ifodd\count0
\def\\{ }\ifx\theshorttitle\relax \thetitle \else\theshorttitle\fi\hfill
{\lnum\number\count0}\else\def\\{ and }{\lnum\number\count0}
\hfill\ifx\theshortauthors\relax 
\theauthors\else\theshortauthors\fi\fi\fi}}\def\@evenhead{\@oddhead}
\def\@oddfoot{\small\lfoot\ifnum\count0=\startpage\copyright\ \gtp\hfill\else
\gt, Volume \thevolumenumber\ (\thevolumeyear)\hfill\fi}
\def\@evenfoot{\@oddfoot}
\makeatother
\fi

\newwrite\gtoutfile
\long\gdef\makeheadfile{  
{\def\\{, }\def\s{ }
\immediate\openout\gtoutfile head.xxx
\immediate\write\gtoutfile{To: math@arxiv.org}
\immediate\write\gtoutfile{Subject: put or rep NNNNN:pppp}
\immediate\write\gtoutfile{--text follows this line--}
\immediate\write\gtoutfile{Proxy-for: \ifx\theasciiauthors\relax
\theauthors\else\theasciiauthors\fi\s<\ifx\theasciiemail\relax\theemail\else\theasciiemail\fi>}
\immediate\write\gtoutfile{\noexpand\\}
\immediate\write\gtoutfile{Authors: \ifx\theasciiauthors\relax
\theauthors\else\theasciiauthors\fi}
\immediate\write\gtoutfile{Title: \ifx\theasciititle\relax
\thetitle\else\theasciititle\fi}
\immediate\write\gtoutfile{Subj-class: GT or SG or MG etc}
\immediate\write\gtoutfile{MSC-class: \theprimaryclass\ifx\thesecondaryclass\relax\else, \thesecondaryclass\fi}
\immediate\write\gtoutfile{Journal-ref: Geom. Topol. \thevolumenumber
(\thevolumeyear) \startpage-\finishpage}
\immediate\write\gtoutfile{Comments: Published by Geometry and Topology at}
\immediate\write\gtoutfile{\s\s http://www.maths.warwick.ac.uk/gt/GTVol\thevolumenumber/paper\thepapernumber.abs.html}
\immediate\write\gtoutfile{\noexpand\\}
\immediate\write\gtoutfile{}
\ifx\theasciiabstract\relax
\immediate\write\gtoutfile{\theabstract}\else
\immediate\write\gtoutfile{\theasciiabstract}\fi
\immediate\write\gtoutfile{}
\immediate\write\gtoutfile{\noexpand\\}
\immediate\write\gtoutfile{}
\immediate\closeout\gtoutfile}}  

\def\maketitlepage{\maketitlep\makeheadfile}
\let\maketitle\maketitlepage

\lognumber{106}
\volumenumber{5}\papernumber{19}\volumeyear{2001}
\pagenumbers{579}{608}
\received{7 January 2000}
\revised{13 June 2000}
\accepted{4 June 2001}
\published{18 June 2001}
\proposed{Gang Tian}
\seconded{Ronald Stern, Ronald Fintushel}

\let\cal\mathcal
\usepackage{amssymb}
\usepackage{euscript}
\usepackage{epsfig}

\newcommand{\pp}{\mathbb{P}}
\newcommand{\qq}{\mathbb{Q}}
\newcommand{\cc}{\mathbb{C}}
\newcommand{\rr}{\mathbb{R}}

\newcommand{\zz}{\mathbb{Z}}
\newcommand{\sss}{\mathbb{S}}

\newcommand{\nn}{\mathbb{N}}
\newcommand{\ff}{\mathbb{F}}

\DeclareFontFamily{OT1}{rsfs}{}
\DeclareFontShape{OT1}{rsfs}{n}{it}{<-> rsfs10}{}
\DeclareMathAlphabet{\curly}{OT1}{rsfs}{n}{it}

\newcommand{\mgbar}{\overline{M}_g}

\newtheoremstyle{Plain}{14pt plus6.3pt minus6.3pt}{7.4pt plus3pt minus3pt}%
{\sl}{}{\bf}{}{0.75em}{\thmnumber{\rm(#2)\bf}\thmname{\hskip0.5em#1}%
\thmnote{\rm\stdspace(#3)}}
\newtheoremstyle{Definition}{14pt plus6.3pt minus6.3pt}%
{7.4pt plus3pt minus3pt}{\rm}{}{\bf}{}{0.75em}{\thmnumber{\rm(#2)\bf}%
\thmname{\hskip0.5em#1}\thmnote{\sl\stdspace#3}}

\theoremstyle{Plain}  
\newtheorem{Thm}{Theorem}[section]
\newtheorem{Prop}[Thm]{Proposition}
\newtheorem{Lem}[Thm]{Lemma}
\newtheorem{Cor}[Thm]{Corollory}

\theoremstyle{Definition}

\newtheorem{Defn}[Thm]{Definition}
\newtheorem{Conj}[Thm]{Conjecture}
\newtheorem{Example}[Thm]{Example}
\newtheorem{Rmk}[Thm]{Remark}

\newtheorem{Question}[Thm]{Question}

\newenvironment{Pf}{\proof}{\endproof}

\newenvironment{Eqn}{\refstepcounter{Thm} $$} {\leqno{\rm(\theThm\rm)} $$} 

\begin{document}

\title{Lefschetz pencils and divisors in moduli space}
\author{Ivan Smith}
\address{New College, Oxford\\OX1 3BN, UK}
\email{smithi@maths.ox.ac.uk}
\begin{abstract}We study Lefschetz pencils on symplectic
  four-manifolds via the associated spheres in the moduli spaces of curves,
  and in particular their intersections with certain natural
  divisors.  An invariant defined from such intersection numbers can
  distinguish manifolds with torsion first Chern class.  We prove
  that pencils of large degree always give spheres which behave
  `homologically' like rational curves;  contrastingly, we give
  the first constructive example of a symplectic non-holomorphic Lefschetz
  pencil.  We also prove that only finitely many values of signature
  or Euler characteristic are realised by manifolds admitting
  Lefschetz pencils of genus two curves.
\end{abstract}

\maketitlepage


\section{Statements of results}

The following section shall set the investigations of this paper into
a wider context, but we record the main results here for the
convenience of the reader.  Recall from \cite{ivanhodge} that a
symplectic four-manifold gives rise to a sequence of spheres $\sss^2
_{X,k} \rightarrow \overline{M}_{g(k)}$ indexed by an integer $k$.
The \emph{covering sequence} of $X$ is the sequence of rational intersection
numbers $(\sss^2 _k \cdot D_k)_{k \in 2\zz}$.  Here $D_k^1 \subset
\overline{M}_{2k-1}$  is a
rational multiple of the divisor of curves which are $k$--fold covers of
$\pp^1$, lifted to
$D_k \subset \overline{M}_{2k-1,h(k)}$ and translated by a term
$\psi_{h(k)}$. The precise definition is given in
(\ref{coveringsequence}).  Suppose that $X$ is a symplectic manifold
whose minimal 
model is not rational or ruled.

\begin{Thm} \label{coneiszero} If the covering sequence of $X$ is
  bounded above then $2c_1 
  (X)$ $ = 0$.  It vanishes identically iff $2c_1 (X)
  = 0 = c_2 (X)$. 
\end{Thm}

\noindent This is a disguised version of results of Taubes and others
  on the sign 
of $K_X \cdot \omega$ for most symplectic manifolds.  The bulk of
the proof involves showing that the relevant Lefschetz pencils
\emph{contain no reducible fibres}, which is a vanishing result for
certain intersection numbers of the sphere in moduli space.  

\begin{Cor}  \label{welose}  For any symplectic manifold $X$ which is
  not rational 
  ruled, and pencils of high
  degree $k$ on $X$, the sphere $\sss^2 _{X,k}$ meets all known effective
  divisors algebraically positively.
\end{Cor}

\noindent In this sense the spheres are ``homologically rational''.  The
  results for small rather than asymptotically large degree $k$ are
  more satisfactory. 

\begin{Thm}  \label{foundone} There is a symplectic genus three
  Lefschetz pencil which is not holomorphic.
\end{Thm}

\noindent In fact the positive relation
  (\ref{Fullerrelation}) lifts to the 
  once marked mapping class group and defines such a pencil.  This is
  the first existence proof for symplectic non-holomorphic pencils
  independent of Donaldson's theorem.  All previous (explicit) examples
  of symplectic non-K\"ahler
  Lefschetz fibrations arose from fibre sum operations and admitted no
  sections of square $(-1)$.  We do not
  establish whether the total space is in fact K\"ahler; it is
  homeomorphic to a complex surface of general type.  Moreover we
  have:

\begin{Thm} \label{genustwofiniteness} Only finitely many pairs $(c_1
  ^2, c_2)$ are realised by 
  the total spaces of genus two Lefschetz pencils.
\end{Thm}

\noindent This is false for Lefschetz fibrations of genus two.  It
  would be a non-trivial consequence of the
symplectic isotopy conjecture of Siebert and Tian, for which it
therefore adduces additional evidence.  It remains a very interesting
question to investigate corresponding finiteness statements for the
geography of manifolds with pencils of higher genus.

\vspace{0.3cm}

\noindent \textbf{Acknowledgements}\qua  Thanks to Denis Auroux for 
comments on an earlier draft of this paper  and to Terry Fuller for
permission to reproduce one of his examples.  This material formed the
basis of a talk at the 2000 Gokova Conference on Geometry and Topology
-- I am grateful to the organisers and participants for a very enjoyable
conference.  Special thanks to Justin Sawon for help taming tortoises.


\section{Introduction}

A fundamental problem in four-dimensional topology is to
establish the relationship between arbitrary smooth four-manifolds and
symplectic four-manifolds on the one hand, and between symplectic
four-manifolds and K\"ahler surfaces on the other.  To this end, the
following questions are natural:

\begin{Question} \label{Question}  (1)\qua   Given an almost complex four-manifold
  $X$ and $\alpha \in H^2 (X; \rr)$ of
  positive square, does $\alpha$ contain symplectic forms inducing the
  given almost complex structure?  

(2)\qua   If the symplectic manifold $X$
  is homeomorphic to a K\"ahler 
  surface, is the symplectic form isotopic to a K\"ahler form?  

(3)\qua  If
  $X$ is K\"ahler and $C \subset X$ is a symplectic submanifold
  realising a homology class with smooth complex representatives, is
  $C$ itself isotopic to a complex curve?
\end{Question}

\noindent  Negative examples for each of these questions are known,
  largely following work in gauge theory.  However, until 
recently no other techniques had yielded comparable progress.  (The
appeal of insight from other arenas is clear, given the continuing
mystery of the analogous questions in higher dimensions.)  In the
papers \cite{SKDIS} and \cite{ivanSDforSS}, Simon Donaldson and the author
apply the machinery of
Lefschetz pencils on symplectic manifolds to construct symplectic
submanifolds, reproving some results of Taubes.  In particular this
leads again to negative examples for the first question.  This paper
describes one 
unsuccessful attempt to use Lefschetz pencils to answer the second and
third questions.  Despite the lack of success, we believe the methods
are of interest.

\noindent The paper is essentially a sequel to \cite{ivanhodge} but to be
self-contained  we
begin with a few background notions.

\begin{Defn}
A \emph{Lefschetz pencil} $f\co X \dashrightarrow \sss^2$ on a four-manifold $X$
comprises a map $f$ from the complement of finitely many points $p_i$
in $X$ to the two-sphere, with finitely many critical
points $q_j$, all in distinct fibres, such that 
\begin{enumerate}
\item $f$ is locally quadratic
near all its critical points $q_j$:  there are local complex co-ordinates
with respect to which the map takes the form $(z_1, z_2) \mapsto z_1
z_2$;
\item $f$ is locally Hopf near all its base-points $p_i$:  there are
  local complex co-ordinates with respect to which the map takes the
  form $(z_1, z_2) \mapsto z_1 / z_2$.
\end{enumerate}
\noindent Moreover all the local complex co-ordinates must agree with
fixed global orientations.
\end{Defn}

\noindent It follows that the
four-manifold is symplectic (we will say more about the cohomology
class of the symplectic form below, cf \ref{notalwayseasy}).  In
projective geometry the appearance of
Lefschetz pencils is classical;  a
generic pencil of divisors on a complex surface gives rise to such a
structure, and blowing up the base-points yields a ``Lefschetz
fibration''.  More recently the inspirational work of Donaldson
\cite{Donaldson:pencils} has
extended the techniques and the descriptions to the symplectic
category:

\begin{Thm}[Donaldson] \label{pencilsexist}
Let $(X, \omega)$ be an integral symplectic manifold and let
$\mathcal{L}_{\omega}$ denote the line bundle with first Chern class
$[\omega]$.  
\begin{enumerate}
\item For $k \gg 0$ there exist pairs $(s_1, s_2)$ of
approximately holomorphic sections of $\mathcal{L}_{\omega} ^{\otimes
  k}$ such that 
the map $X \backslash \{s_1 = s_2 = 0 \} \rightarrow \pp^1$ may be
perturbed to define a Lefschetz pencil.  
\item The fibres of the pencil are symplectic submanifolds (away from
  the finite set of critical points) representing the Poincar\'e dual
  of $k[\omega]$ in $H_2 (X, \zz)$.
\item Once $k$ is
sufficiently large, the pencils obtained in this way are canonical up
to isotopy.
\end{enumerate}
\end{Thm}

\noindent If we blow up the base-points $p_i$ of a pencil, the map $f$
extends to the total space and we obtain a
Lefschetz fibration.  The generic
fibre of the 
fibration is a smooth two-manifold $\Sigma_g$ of some fixed genus $g$;
in this paper $g \geq 2$ unless stated otherwise. 

\begin{Rmk}[On notation]
We will use
the term ``fibre'' to refer to the gen\-eric hypersurface in a pencil,
as well as the preimage of a point in a fibration.
Note that we reserve the term
``pencil'' to refer to the four-manifold before blowing up, and so the
fibres in a pencil have strictly positive square.
\end{Rmk}

\noindent The
topology of the four-manifold $X$ is encoded in a \emph{positive
  relation}; this is a word 
$\langle \prod \delta_i = 1 \rangle$ in
positive Dehn twists in the mapping class group
$\Gamma_g$ which
determines the monodromy representation and the diffeomorphism
type of the fibration.  Here $\delta_i$ is the positive Dehn twist
about some fixed (isotopy class of) embedded curve $C_i \subset
\Sigma_g$;  occasionally we shall use $\delta$ for the curve as well
as the twist diffeomorphism.  Such curves are called \emph{vanishing
  cycles} and generate a subgroup $V < \pi_1 (\Sigma_g)$, which we
  always assume is non-empty.  The fundamental group of the 
four-manifold $X$ is given by 
$\pi_1 (\Sigma_g) / V$.  All of our vanishing cycles $C$ will be
  homotopically essential and hence the fibres will contain no spherical
  components.  Given this, if we choose a metric on $X$ the
  smooth fibres 
become Riemann surfaces and the critical fibres stable Riemann
surfaces, and we induce a map $\phi_f\co  \sss^2 \rightarrow \mgbar$ with
image the Deligne--Mumford moduli space of stable curves.  Recall this
moduli space is given by adjoining certain divisors of stable curves
$\Delta_i$ ($0 \leq i \leq [g/2]$) to the moduli space of smooth curves $M_g$.
The $\Delta_i$ form
the irreducible components of a connected divisor $\Delta = \cup \Delta_i$ and
the generic curve in $\Delta_i$ has one component of genus $i$ and one of
genus $g-i$ if $i>0$, and is irreducible if $i=0$.  The fibration
  being symplectic
means that the sphere $\sss^2 _X$ has \emph{locally positive
  intersections} with 
the various divisors $\Delta_i$; the restriction to nodal
  singularities means the intersections are \emph{transverse}.  In this
  paper we shall refer to a fibration $f$ on $X$ giving rise to a
  sphere $\sss^2 _X$ and will suppress the choice of metric; changing
  the metric changes the sphere by an isotopy which always preserves
  the geometric intersection number with the $\Delta_i$ (an
  ``admissible isotopy'').  All
  results depend only on the admissible isotopy class of the sphere.
  
\begin{Rmk}
In \cite{SKDIS} and \cite{ivanSDforSS}, symplectic submanifolds are
constructed by studying
the Gromov invariants of fibrations of symmetric products associated
to a Lefschetz pencil.  In \cite{ivaninstanton} Gromov invariants for
associated families of moduli spaces of stable bundles are related to
instantons on the four-manifold.  However, the spheres $\sss^2
\rightarrow \mgbar$ in this paper are \emph{not necessarily symplectic
  or pseudoholomorphic}.  It remains an interesting question to study
the quantum cohomology of moduli spaces of curves.
\end{Rmk}

\noindent The geometric classes $\Delta_i$ define
  elements of $H^2 (\mgbar; \zz)$; there is an important algebraic
  class which is
  not dual to any distinguished divisor.  This is the Hodge class
  $\lambda$, which is the first Chern class of the relative dualising
  sheaf of the universal curve $\mathcal{C}_g \rightarrow M_g$.
  Although the universal curve does not exist, the Chern class makes
  sense.  There are
  additional natural classes when we look at moduli spaces of curves
  with (ordered) marked points.  Let $\psi_i$ be the
first Chern class of the line bundle on $\overline{M}_{g,h}$ with
fibre $T^*_{p_i}$ at $(\Sigma; p_1, \ldots, p_h)$.  We record
  the following, which is due to Harer \cite{Harer} and
  Arbarello--Cornalba \cite{arabel}:

\begin{Thm}[Harer, Arbarello--Cornalba]
The classes $\Delta_i$ with $0 \leq i \leq [g/2]$ and $\lambda$
rationally generate the second cohomology $H^2 (\mgbar; \qq)$.  For
the pointed moduli space, $H^2 (\overline{M}_{g,h}; \qq)$ is rationally
generated by the pullbacks of these classes under the 
forgetful map and the classes $\{ \psi_i \ | \ 1 \leq i \leq h \}$.
\end{Thm}

\noindent The two descriptions of a four-manifold given by a positive
relation and a sphere in moduli space seem rather far from one
another.  Nonetheless, they may be related, and we shall exploit the
dual descriptions in 
investigating the consequences of the following trivial result:

\begin{Prop} \label{principle}
If $f\co X \rightarrow \sss^2$ gives rise to a sphere $\sss^2 \subset
\mgbar$ and for some effective divisor $D \subset \mgbar$ not containing
$\sss^2$ we have $[\sss^2]\cdot[D] < 0$ then $f$ is not isotopic to a
holomorphic fibration.
\end{Prop}

\noindent This suggests a natural obstruction to the existence of
K\"ahler forms isotopic to given symplectic forms:

\begin{Cor} \label{obstruction}
Let $(X, \omega)$ be a symplectic manifold and suppose the
Lefschetz pencils obtained from asymptotically holomorphic sections of
$\mathcal{L}_{\omega}^{\otimes k}$ give spheres
$\sss^2 _k \subset \overline{M}_{g(k)}$.  If for some sequence of
divisors $D_k \subset \overline{M}_{g(k)}$ we have that $\sss^2 _k
\not \subset D_k$ and $[\sss^2 _k] \cdot D_k < 0$ then $\omega$ is not
deformation equivalent to a K\"ahler form on $X$.
\end{Cor}

\begin{Pf}
According to Donaldson \cite{Donaldson:pencils} the pencils defined by his
construction from approximately holomorphic sections are isotopic to the
pencils provided by algebraic geometry for $k \gg 0$.  These pencils
give rise to rational curves in moduli space, which must meet
positively all
effective divisors in which they are not contained.  Thus the
assumptions indicate that the spheres defined by the given manifold
$X$ are not isotopic to rational curves. By the asymptotic uniqueness
in Donaldson's theorem, this shows that the symplectic structure on
$X$ is not in fact K\"ahler.
\end{Pf}

\noindent It is important to note that these obstructions might well
be computable.  Given one explicit positive relation (with
base-points) defining a pencil, there is a stabilisation process which
obtains pencils of higher degree \cite{AKstabilisation}.  If the
initial pencil is obtained from approximately holomorphic sections, so
are the later ones.  Moreover, there are familiar techniques for
computing cohomology classes of divisors in moduli space and hence the
intersection numbers above, whilst the condition that the sphere lies inside a
particular divisor may have topological consequences which can be
checked independently.  We shall use such an argument in the proof of
(\ref{foundone}) later in the paper.  Nonetheless, our main
observation -- contained in the theorems (\ref{coneiszero}) and
(\ref{welose}) -- amounts to the 
triviality of these obstructions.  In principle it remains possible
that the obstructions could detect symplectic manifolds with $2c_1
= 0$ which were not K\"ahler, but this seems implausible.

\noindent In the face of these results, we
shall turn from studying the asymptotic intersection behaviours to
concentrating on rather particular divisors at small genus;  in this
framework we shall deduce the results (\ref{foundone}) and
(\ref{genustwofiniteness}).  The first instance involves 
the hyperelliptic divisor inside the moduli space of genus
three curves.  In particular, we show that in certain circumstances
one can detect \emph{from a positive relation} the non-holomorphicity of a
fibration, even when the total space is homotopy K\"ahler.  The
finiteness result will rely on
working with pairs comprising a Lefschetz fibration and a
distinguished section, and with a certain Weierstrass divisor in the
moduli space of pointed genus two curves.  We leave the details until
the relevant discussion in the paper.  We do draw attention to one
general fact however.  The first constructions of symplectic
structures on manifolds with Lefschetz pencils involved blowing up to
the Lefschetz fibration, and blowing down exceptional symplectic
sections again.  This led to an unfortunate asymmetry:  Donaldson's
theorem (\ref{pencilsexist}) gives Lefschetz pencils whose fibres are
dual to a given integral symplectic form, but given a Lefschetz pencil
there may be \emph{no symplectic form} in the integral cohomology
class which is dual to a fibre.

\begin{Lem} \label{notalwayseasy}
There are smooth four-manifolds $X$ with the topological structure of
Lefschetz pencils for which the class $PD[\mathrm{Fibre}] \in H^2 (X;
\zz)$ admits no symplectic representative symplectic on the fibres.
\end{Lem}

\begin{Pf}
Let $X$ be a manifold with a
pencil of curves with one base-point, and for which some member of the
pencil is a reducible curve.  Such can be obtained by blowing up all
but one of the base-points on any pencil containing reducible elements,
for instance the pencil obtained by Matsumoto \cite{Matsumoto}.
Suppose for contradiction that there is a symplectic form in the class
$PD[\mathrm{Fibre}]$ for which the smooth locus of each fibre is a
symplectic submanifold.  Then each component of a reducible fibre has
positive symplectic area, and hence since the symplectic form is
integral $[\omega] \cdot [\mathrm{Fibre}] \geq 2$.  But we have
assumed there is a unique base-point.
\end{Pf}

\noindent This is the only source of such pathologies:  according to a
recent theorem of Gompf \cite{GompfGokova}, if there are base-points on
every component 
of every fibre then there is a symplectic form in the distinguished
integral class.  It follows that the non-K\"ahler pencil that we
construct will indeed determine a canonical isotopy class and not
deformation equivalence class of symplectic form on the
four-manifold.

\noindent  We note in closing that there is an analogue of the foundational
(\ref{principle}) which applies to branched coverings:

\medskip

  {\sl Suppose $\overline{\mathcal{N}}_g \rightarrow \mgbar$ is a
  branched covering.  Then a sphere $\sss^2 \subset \mgbar$ lifts to
  the covering space $\overline{\mathcal{N}}_g$ if and only if it
  meets the branch locus everywhere tangentially.}

\medskip

\noindent There is an obvious topological constraint which
must be satisfied if such tangential intersections can arise: the
intersection number between the sphere and the branch divisor must be
\emph{even}.  This leads to a circle of ideas somewhat similar to that
developed in this paper.


\section{Asymptotic intersections}
  
In this section we shall prove the result (\ref{coneiszero}) and
explain the shape of the statement (\ref{welose}).  
There is a well-known conjecture on the geometry of the moduli spaces
$\mgbar$.  Recall Harer's theorem from above;  the following is taken
from \cite{moduliofcurves}.

\begin{Conj}[``Slope Conjecture'', Harris--Morrison\/]
Let $a\lambda - b\Delta_0$ be\break represented by an effective divisor
in $\mgbar$.  The ratio $a/b$ is minimised by  Brill--Noether
divisors.
\end{Conj}

\noindent This special property of Brill--Noether
divisors\footnote{That is, divisors of curves which have a linear
  system they shouldn't have.} motivates
the definition of a particular symplectic invariant.  By the
adjunction formula, the pencils of even degree $2k$ on any symplectic
manifold give rise to families of curves of odd genus.  In these
moduli spaces we have divisors
$D^1_{(g+1)/2} \subset \mgbar$ comprising the
codimension one components of the closure in $\mgbar$ of the locus of
curves in $M_g$ which admit a $g^1_{(g+1)/2}$.  (Recall that a $g^r_d$
  is a linear system of dimension $r$ and degree $d$.)  Since a
linear system of dimension one and degree $d$ is just a $d$--fold
covering over $\pp^1$, the hyperelliptic divisor is precisely $D^1_2
\subset \overline{M}_3$.  For even $g$ there are analogous ``Petri''
divisors, but we will avoid the notational complication they
introduce. 

\begin{Thm}[Harris--Mumford] \label{Brill-Noether}
Let $g$ be odd and $k=(g+1)/2$.  The cohomology class of $D^1 _k$ is
given by
$$D^1 _k \ = \ c_k \big( (g+3) \lambda - \frac{g+1}{6} \Delta_0 -
\sum_{i=0}^{[g/2]} i(g-i) \Delta_i \big).$$
\noindent Here $c_k = 3(2k-4)! / ( k! (k-2)! )$ is positive and rational.
\end{Thm}

\noindent It shall be important for us to normalise the Brill--Noether
divisors by dividing by these rational constants $c_k$.  A Lefschetz
pencil with  $h$ base-points gives rise to a
sphere in a moduli space $\overline{M}_{g,h}$ of  curves with $h$
ordered marked points.  We need to keep track of the marked points in
order to utilise all the geometry of the system; we do this by
translating the divisor $D^1 _k$ when we lift to $\overline{M}_{g,h}$.

\begin{Defn} \label{coveringsequence}
Let $(X, \omega)$ be an integral symplectic manifold.  The
\emph{covering sequence} of $(X, \omega)$ is the sequence of
intersection numbers $([\sss^2_k] \cdot D_k)_{k \in 2\nn}$ between
the spheres defined by pencils of curves dual to $k\omega$ on $X$ and
the divisors
$$D_k \ = \ \frac{1}{c_k} D^1 _k - \sum_{j=1}^{k^2 \omega^2} \psi_j.$$
\noindent We set the intersection
number to be zero by convention if there is no pencil of curves in the
relevant homology class, or if $g=g(k)$ is not odd.
\end{Defn}

\noindent If two integral symplectic manifolds $(X, \omega_X)$ and $(Y,
\omega_Y)$ are symplectomorphic then it follows from the theorem of
Donaldson (\ref{pencilsexist}) that their covering sequences co-incide
for sufficiently large $k$.  The values of the covering sequence compare the
numbers of exceptional (for instance hyperelliptic) fibres in the
Lefschetz fibration, normalised by the universal constants $c_k$, with
the number of exceptional sections of the fibration.  This
is quite a natural comparison:  at sufficiently large $k$ we know the
pencils of approximately holomorphic sections can be extended to nets
(two dimensional linear systems) of sections \cite{Auroux}.  A net
gives a branched cover of $X$ over $\cc \pp^2$, or equivalently a
branched cover of the total space of a Lefschetz fibration over the
first Hirzebruch surface $\ff_1$ in which all the exceptional sections of the
Lefschetz fibration map to the unique holomorphic section of square
$(-1)$ in $\ff_1$.  Hence for nets to exist, we know that the fibres
of the Lefschetz pencil must admit branched coverings over $\pp^1$ of
the explicit degree given by the number of base-points.  Thus we
are comparing this degree valid for all the fibres to the minimal
degree realised by some of them.

\noindent To establish a link between four-manifolds and intersection
numbers, we shall need to understand the values of the generators of
$H^2$ on a sphere $[\sss^2]$.  For the $\Delta_i$ it is
straightforward -- we count the numbers of singular fibres of different
kinds.  The main result of \cite{ivanhodge} solves the problem for the
Hodge class, which is related to the signature of the associated
Lefschetz fibration:  
\begin{Eqn} \label{signatureformula}
\sigma(X) = 4 \langle \lambda, [\sss^2 _X]
\rangle - \sum \langle \Delta_i, [\sss^2 _X] \rangle.
\end{Eqn}
\noindent We have an easy lemma:

\begin{Lem} \label{sectionsquare}
Let $X \rightarrow \sss^2$ be a Lefschetz fibration with a
distinguished section $s$.  Then the value $\langle \psi,
[\sss^2_{X,s}] \rangle$ of $\psi \in H^2 (\overline{M}_{g,1}; \zz)$ on
the associated sphere is given by the self-intersection $s \cdot s$.
\end{Lem}

\noindent Lastly recall the adjunction formula;  for a pencil of curves on $X$
dual to $k[\omega]$ we have
\begin{Eqn} \label{adjunction}
2g-2 \ = \ K_X \cdot k[\omega] + k^2 [\omega]^2.
\end{Eqn}
\noindent From these pieces of information, and the result
(\ref{Brill-Noether}), we can now derive all the values of the
intersections $\sss^2 _k \cdot D_k$ from a combinatorial description
of a Lefschetz fibration. However, the topological type of the
manifold does not suffice, largely because it is not clear how to
determine the individual values $\Delta_i$ rather than their sum.  In
fact, if $i > 0$ we can expect the value to be zero.  The motivation
comes from complex geometry:

\begin{Prop} \label{Hodgeindextheorem}
Suppose the K\"ahler surface $X$ has a Lefschetz pencil of curves in a
homology class $[C]$.  If the pencil contains a reducible curve, then
either $[C]$ is
indivisible in $H_2 (X, \zz)$ or $[C] = 2[D]$ with $[D]$ indivisible and
$[D]^2 = 1$.
\end{Prop}

\begin{Pf}
Recall that the Hodge Index theorem asserts that the
intersection form on $H^{1,1}$ has a unique positive eigenvalue for
any K\"ahler surface.  In particular, given any divisor $D$
for which $D^2 > 0$ and a divisor $D'$ such that $D \cdot D' = 0$ we
can deduce that either $D' = 0$ or $D' \cdot D' < 0$.   We apply this
in the following way.  First note, scaling $k$ if necessary, we have
in the notation  of (\ref{reducibleones}) that $(D_1)
^2 + (D_2) ^2 = k^2 [C] - 2 > 2$, and we may assume without
loss of generality that $(D_1) ^2 > 0$.  Now 
$$D_1 \cdot (D_1 - (D_1)^2 D_2) = 0 \ \Rightarrow \ D_1 = (D_1)^2 D_2
\ \ \mathrm{or} \ \ (D_1 - (D_1)^2 D_2)^2 < 0.$$
\noindent The first case gives an easy contradiction.  If $D_1 =
(D_1)^2 D_2$ then $D_1 \cdot D_2 = (D_1)^2 (D_2)^2 = 1$ which forces
$(D_1)^2 = (D_2) ^2 = 1$ by integrality; but this is impossible, since
$(D_1)^2 + (D_2)^2 > 2$.  So assume instead that the second case
holds, and expanding we find that $(D_1)^2 [(D_1)^2 (D_2)^2 - 1] < 0$
which gives $(D_1)^2 (D_2)^2 < 1$, and hence $(D_2)^2 \leq 0$ since
the first of the two terms is assumed positive.  Now since
$$k \, | \, (D_2 \cdot [k \omega]) \ = \ D_2 (D_1 + D_2)  \ = \ 1+(D_2)^2$$
\noindent we know that we cannot have $(D_2)^2 = 0$ as soon as $k$ is
not equal to $1$.  But now we have a contradiction again: by
assumption, each of the $D_i$ appeared as complex curves in a pencil
and hence has positive symplectic area.  From which
$$[k\omega] \cdot D_2 > 0 \ \Rightarrow \ 1 + (D_2)^2 > 0$$
\noindent which is absurd if also $(D_2)^2 \leq -1$.  Thus as soon as
$k$ is large enough that $k^2 [C]^2 > 4$ and $k>1$, no Lefschetz
pencil of curves in the class $[C]$ on a K\"ahler surface can contain reducible
elements.
\end{Pf}

\noindent There are analogues of this result for general symplectic
pencils.  With the stabilisation procedure for pencils in mind, we
shall concentrate on pencils of even degree.  Then there is a trivial
argument for certain classes of
four-manifold, even without the K\"ahler assumption:

\begin{Prop}
Suppose the symplectic manifold $X$ has even intersection
form, or that $K_X$ is two-divisible in cohomology.  Then a
Lefschetz pencil dual to $2k [C]$, for any $k \in \zz_+$ and $[C] \in
H_2 (X; \zz)$, contains no reducible fibres.
\end{Prop}

\begin{Pf}
To obtain
a reducible curve in a pencil of curves dual to $k[C]$ precisely involves
writing
\begin{Eqn} \label{reducibleones}
k[C] = D_1 + D_2 \ \in H^2 (X, \zz) \qquad \mathrm{where} \qquad D_1 
\cdot D_2 = 1.
\end{Eqn}
\noindent But this is impossible under the additional assumptions on $X,k$; for
if $D_1 = r[C] + \alpha$ and $D_2 = (k-r)[C] - \alpha$ then 
$$D_1 \cdot D_2 \ = \ r(k-r) [C]^2 + (k-2r)[C] \cdot \alpha - \alpha^2.$$
\noindent But if $k$ is even and the intersection form is even, all three terms
are divisible by two and hence the LHS cannot equal $1$.  Using another trick,
we see that if $K_X$ is even (ie, can be written as $2\kappa$ for
some cohomology class $\kappa$) then by adjunction $(D_1)^2$ must be
even:
$$2g_{D_1} - 2 \ = \ K_X \cdot D_1 + (D_1)^2$$
\noindent but on the other hand, for even $k$, we have
$$D_1 \cdot (k [C]) \ = \ D_1 \cdot (D_1 + D_2)  \ = \ (D_1)^2 + 1$$
\noindent forcing $(D_1)^2$ to be odd, a contradiction.  
\end{Pf}

\noindent Nonetheless, to obtain the general case requires more work.
The following can be regarded as a symplectic shadow of the Hodge
Index theorem.  Note that it is in fact a vanishing result;  it
implies that the intersection numbers $[\sss^2 _k] \cdot \Delta_i = 0$
whenever $k$ is large and even, and $i > 0$.  A detailed algebraic
treatment of stabilisation is now available in \cite{AKstabilisation}.

\begin{Thm} \label{vanishing}
Every symplectic four-manifold admits Lefschetz pencils\break composed of
only irreducible curves.  Indeed the pencils arising from
the stabilisation $k \mapsto 2k$ procedure always have this form.
\end{Thm}

\begin{Pf}   The central observation we need
  is due to Donaldson;  it is easier to pass from Lefschetz
  pencils representing $k\omega$ to ones representing $2k \omega$ than
  to ones representing $(k+1) \omega$.  Suppose then we have a
  symplectic four-manifold $X$ and a pencil of sections of a line
  bundle $L$ with $c_1 (L) = k\omega$ (normalised appropriately).  Let
  the pencil of sections $\{ s_1 + \lambda s_2 \}_{\lambda \in \pp^1}$
  be generated by two smooth elements $s_1, s_2$.  We consider the
  pencil of reducible nodal sections $\{ s_1^2 + \lambda s_1 s_2
  \}_{\lambda \in \pp^1}$.  Thus we ``add in'' the zero-set of $s_1$
  to each of the curves in the original pencil.  In complex geometry
  this would correspond to taking a sphere $\pp^1 \subset \pp H^0 (L)$
  with image entirely contained in the discriminant locus of singular
  curves.  Such a sphere could be perturbed to have isolated
  transverse intersections with the discriminant, corresponding to a
  deformation to the Lefschetz situation; the above remarks would then
  apply.  We mimic this by (the first rather easy steps of) an
  analysis of the relevant deformation.  Note that the singular
  sections still satisfy the approximate holomorphicity and
  $C$--bounded constraints of \cite{Donaldson:pencils} for suitable constants,
  and it is transversality that we must achieve.

\noindent Away from $(\lambda=0)$ in the $\lambda$--plane, we have a
family of nodal curves parametr\-ised by a disc.  There is a small
$C^2$--deformation of this family, analogous to smoothing a single
nodal curve into the self-connect sum of its normalisation. Thus a
perturbation of the form 
$$\{ (s_1^2 + \epsilon) + \lambda (s_1 s_2 +
\delta)\}_{|\lambda| \gg |\epsilon| + |\delta|}$$
\noindent smooths each of the curves $s_1^2, s_1 s_2$ and for generic
smoothing sections $\epsilon, \delta$ will give at least away from
$\lambda = 0$ a Lefschetz family by the arguments of
\cite{Donaldson:pencils}.  But 
by continuity we can see the
critical members of this family:  they arise (again for large
$|\lambda|$) in one-to-one correspondence with the critical
fibres of the monodromy of
the original pencil $\{ s_1 + \lambda s_2\}$ on $X$.

\noindent It follows that the new monodromy near
infinity in the $\lambda$--plane gives a copy of the old monodromy in a
model where we identify each original fibre $(s_2 =
\mathrm{constant})$ minus small discs
with its image in a smoothed fibre.  The product of the original Dehn
twist monodromies is no longer the identity but equal to a word in the
(commuting) Dehn twists about the necks that link a copy of $s_2$ with
a copy of $s_1$ once the finitely many nodal intersections are
smoothed.  Thus we must determine an integer $n_i$ for each of these
nodes, and the monodromy about a large circle around $0$ in the
$\lambda$--plane is given by $\prod \delta_i ^{n_i}$ where $i$ runs
over the base-points of the original pencil.  By symmetry each of the
$n_i$ will be equal, since this is a local question: we are taking a
section of a line bundle over $\pp^1$ with fibre the tensor product of
the tangent directions at the node.  (This can be identified with the normal
bundle to the 
divisor of stable nodal curves in the moduli space $\overline{M}_g$.)
Since the tangent direction at
the fixed curve is constant, the line bundle is the normal bundle to
the section defined by the base-point.
Since the exceptional sections have normal bundle $\mathcal{O}(-1)$ we
deduce that $n_i= -1 \ \forall i$.

\noindent The situation is now as depicted in Figure
\ref{lambdaplane}.  The remaining critical fibres in the $2k$--pencil
come from deforming the multiple fibre $s_1 ^2 = 0$ to a smooth
locus.  This is again a local consideration and we can proceed in
various ways.  If $b_+ (X) > 1$ we can choose an
\begin{figure}[ht!]
\vspace{0.5cm}
\begin{center}
\begin{picture}(0,0)%
\epsfig{file=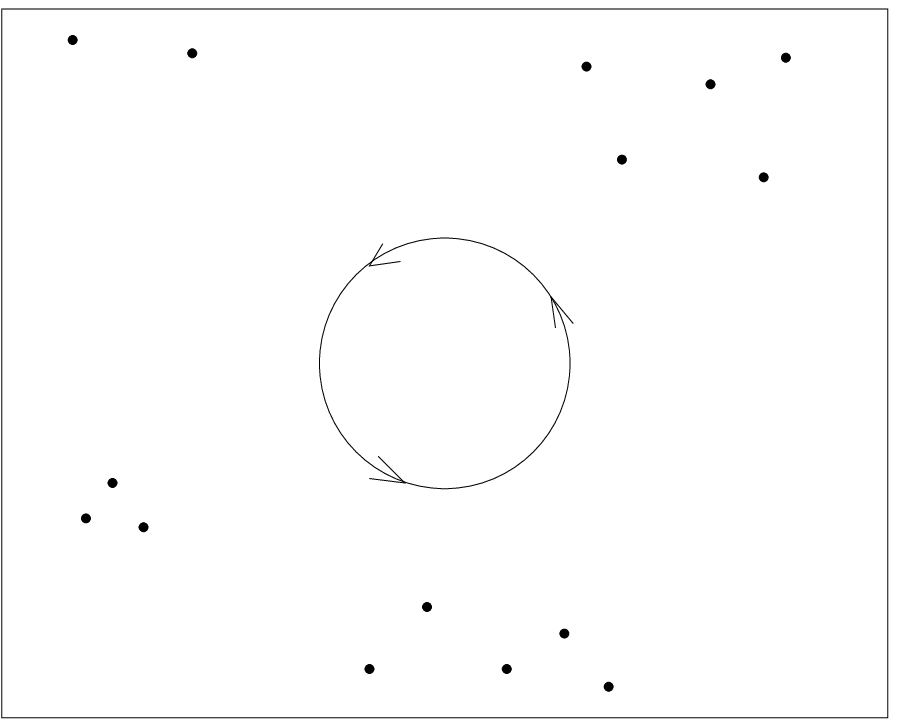}%
\end{picture}%
\setlength{\unitlength}{1776sp}%
\begingroup\makeatletter\ifx\SetFigFont\undefined%
\gdef\SetFigFont#1#2#3#4#5{%
  \reset@font\fontsize{#1}{#2pt}%
  \fontfamily{#3}\fontseries{#4}\fontshape{#5}%
  \selectfont}%
\fi\endgroup%
\begin{picture}(9024,7224)(889,-7273)
\put(4681,-3796){\makebox(0,0)[lb]{\smash{\SetFigFont{5}{6.0}{\rmdefault}{\mddefault}{\updefault}
to be Lefschetz
}}}
\put(2881,-2626){\makebox(0,0)[lb]{\smash{\SetFigFont{5}{6.0}{\rmdefault}{\mddefault}{\updefault}
Circle shows
}}}
\put(7111,-3616){\makebox(0,0)[lb]{\smash{\SetFigFont{5}{6.0}{\rmdefault}{\mddefault}{\updefault}
Monodromy round 
}}}
\put(7426,-3931){\makebox(0,0)[lb]{\smash{\SetFigFont{5}{6.0}{\rmdefault}{\mddefault}{\updefault}
circle: $\prod_i \delta_i ^{-1}$
}}}
\put(1171,-916){\makebox(0,0)[lb]{\smash{\SetFigFont{5}{6.0}{\rmdefault}{\mddefault}{\updefault}
Critical values of initial
}}}
\put(1351,-1141){\makebox(0,0)[lb]{\smash{\SetFigFont{5}{6.0}{\rmdefault}{\mddefault}{\updefault}
monodromy (black dots)
}}}
\put(2386,-2896){\makebox(0,0)[lb]{\smash{\SetFigFont{5}{6.0}{\rmdefault}{\mddefault}{\updefault}
$|\lambda| = |\epsilon| + |\delta|$
}}}
\put(4636,-691){\makebox(0,0)[lb]{\smash{\SetFigFont{6}{7.2}{\rmdefault}{\mddefault}{\updefault}
\textbf{THE $\lambda$-PLANE}
}}}
\put(4681,-3481){\makebox(0,0)[lb]{\smash{\SetFigFont{5}{6.0}{\rmdefault}{\mddefault}{\updefault}
Deform $s_1 ^2 = 0$
}}}
\end{picture}
\caption{Stabilisation viewed in the $\lambda$--plane
  \label{lambdaplane}}
\end{center}
\end{figure}
integrable complex structure in any sufficiently small neighbourhood
of the multiple fibre.  This involves identifying to diffeomorphism
the model with 
a neighbourhood of a multiple fibre in a pencil of singular curves on a
K\"ahler surface -- note that the local diffeomorphism type depends
only on the genus of the curve and the number of base-points of the
original pencil.  On the
other hand, for the K\"ahler situation
reducible curves of Lefschetz type cannot occur for topological
reasons, as described above.  Thus we need to know
that there is always a K\"ahler model available.  If the four-manifold
has $b_+ > 1$ this follows from work of Taubes \cite{Taubes}: we know
that the number of base-points
is bounded above by $2g-2$ (by the positivity of $K \cdot \omega$, cf
\ref{likesymplectic}).  But then take a genus
$g$ pencil on a holomorphic $K3$ surface; this has $2g-2$ base-points,
so blowing these up successively gives pencils of curves with any
strictly positive intermediate number of base-points.  Moreover the
singular holomorphic pencils do have Lefschetz smoothings, by
Riemann--Roch;  in any high dimensional linear system of curves with
some smooth members, the generic pencil is Lefschetz by Bertini's
theorem and complex Morse theory.  For large $k$ we can assume that
our divisors are very ample and smooth curves do exist.

\noindent To avoid the Seiberg--Witten theory, or when $b_+ = 1$, we
can study the local degeneration directly.  For instance, fix a family 
of smooth curves with an arbitrary number $N$ of sections over a disc.
Gluing on a copy of the central fibre along the sections gives a disc
in a moduli space $\overline{M}_{g,N}$ which is contained entirely in
the stable locus.  By projectivity of this stable divisor, there is
some complex curve in the stable locus which contains an isotopic copy
of some sufficiently small sub-disc; then a holomorphic perturbation
of this complex curve 
will be transverse to the stable locus and will also give a model for
the relevant local degeneration.  In another direction, we could model
the $k \mapsto 2k$ stabilisation by a perturbation 
$$\{s_1 + \lambda s_2 \}_{\lambda \in \pp^1} \ \Rightarrow \ \{(s_1 s_3
+ \epsilon) + \lambda (s_2 s_3 + \delta) \}_{\lambda \in \pp^1}.$$
\noindent Here we obtain the section $s_3$ from Auroux's construction
of nets \cite{Auroux}.  By uniqueness at large enough $k$, this pencil
must be isotopic to that given by the previous stabilisation, and the
distinct local models imply that the smoothing of the double fibre
yields no reducible vanishing cycles.

\noindent The proof is completed with a familiar
trick.  For the 
reducible fibres in the original pencil, both
components are smooth symplectic submanifolds themselves.  It follows
that for each component $D_i$ we have $(\omega) \cdot D_i > 0$ which
means that there are base-points of the original pencil on both
components.  Once we add in a copy of $s_1$ and smooth the nodes,
this precisely means that the separating vanishing cycles no longer
separate the new generic fibre.  Thus the separating vanishing cycles
from the old monodromy become non-separating, and the stabilisation is
elsewhere modelled on a complex situation in which there are no
reducible curves.  The theorem follows.
\end{Pf}

\noindent It may be of independent interest to observe that there is a
stabilisation procedure in which the vanishing cycles at level $k$
form a subset of a natural set of vanishing cycles at level $2k$.
Note that we could not make sense of this without considering the
base-points of the pencil as well as the vanishing cycles of the
fibration.  We can now prove the first theorem (\ref{coneiszero}).

\begin{Thm}
Suppose that $X$ is not rational or ruled.  The covering sequence of
$X$ is bounded above only if $2c_1 (X) = 0$, and vanishes
identically iff $2c_1 (X) = 0, c_2 (X) = 0$.
\end{Thm}

\begin{Pf}
We compute the intersection number $[\sss^2 _k] \cdot D_k$ using the
formulae (\ref{Brill-Noether}, \ref{signatureformula},
\ref{adjunction}) as well  as the
definition (\ref{coveringsequence}) and the vanishing theorem
(\ref{vanishing}).  The classes $\psi_i$ evaluate, by the lemma
(\ref{sectionsquare}), on any sphere arising from a pencil to give
$-1$, and hence we can compute the intersection with $D_k^1 / c_k$
and then adjust by subtracting the number of base-points of the
pencil.  Of course this is just $k^2 \omega^2$.  The calculation gives
$$[\sss^2 _k] \cdot D_k \ = \ \frac{(g+7)(g+1)}{12} K_X \cdot \omega +
\frac{g+3}{12} (c_1 ^2 (X) - c_2 (X)) - \frac{2}{3} c_2 (X). $$
\noindent Here $\omega$ refers to the given integral symplectic form
on $X$ and not some multiple, and we have substituted $g = 2k-1$.
This expression will grow positively as $o(g^2)$ unless $K_X \cdot
\omega \leq 0$.  But we know from Seiberg--Witten theory \cite{LiuSW},
\cite{McDS2} that $K_X \cdot \omega < 0$ only if $X$ is a rational or
ruled surface, which we exclude by assumption, and that $K_X \cdot
\omega = 0$ only if $2 K_X = 0$.  If $b_+ > 1$ then we deduce
that $K_X = 0$ itself.

\noindent If $K_X$ is a torsion class, then the intersection number is
given by a negative multiple of $c_2 (X)$ and hence the covering sequence is
identically zero iff this also vanishes.  The last thing we need to
know is that when $K_X$ is torsion, $c_2 (X) \geq 0$ so that the
sequence is necessarily bounded above in this case.  But if $2c_1 (X)
= 0$ and $c_1 (X) \neq 0$ then we know that $b_+ (X) = 1$, and then
the result follows since $0 = c_1 ^2 = 2c_2 + 3(1-b_-)$.  (It is
likely, but not proven in general, that whenever $c_1 (X) = 0$ then
$c_2 (X) \geq 0$; this would give ``iff'' in the first statement of the
Theorem.)  
\end{Pf} 

\noindent The intersection formula
above was twisted by subtracting the large number of
exceptional sections from the naive intersection $[\sss^2 _k] \cdot
[(1/c_k)D^1_k]$).  Given the slope conjecture, we see that for $K_X
\cdot \omega > 0$ the 
intersection numbers of spheres in moduli space with any effective
divisor grow unbounded to infinity.  Indeed this growth is quadratic
with the degree $k$.  This justifies the corollary
(\ref{welose}), and explains the failure of the obstructions
(\ref{obstruction}) to distinguish symplectic and K\"ahler structures
in four dimensions.

\begin{Rmk} \label{generaltype}
If we have a K\"ahler surface  with $K_X \cdot \omega < 0$ then we
have shown that the
spheres defined by pencils of large
enough degree will necessarily be contained in all of the
Brill--Noether divisors.  This however is not so surprising.
Recall that for $g>23$  the moduli space $\mgbar$ is
known to be of general type and not unirational (as the first few
moduli spaces are, when $g \leq 11$).  It follows that through the
general point of $\mgbar$ there is no rational curve, and
\emph{all} the spheres defined by holomorphic fibrations lie
inside distinguished subvarieties of special curves (with no condition on
$K \cdot \omega$ of the underlying surface). 
\end{Rmk}

\noindent In \cite{ivanhodge} we proved that the symplectic area of
the spheres $\sss^2 _k$ is positive for any symplectic manifold $X$
and pencils of any degree $k$.  This follows from the statement that
the evaluation $\langle \lambda, [\sss^2] \rangle$ is strictly
positive.  In the light of the above, it makes sense to ask if the
principle (\ref{principle}) can \emph{ever} be applied;  can any
sphere meet any effective divisor negatively?  Certainly local
negative intersections must arise, for spheres which are not isotopic
to rational curves, but it is not clear that these can ever contribute
sufficiently to give a negative algebraic intersection.  In the next
section we shall provide an explicit example to show that (\ref{principle})
is not entirely vacuous.


\section{Symplectic non-K\"ahler Lefschetz pencils}

Every K\"ahler surface admits a holomorphic Lefschetz pencil.  Is
every Lefschetz pencil on a smooth four-manifold in fact holomorphic
for some K\"ahler structure?  The answer is clearly no, for
Donaldson's theorem
(\ref{pencilsexist}) provides pencils on manifolds not homotopy
equivalent to complex surfaces.  In this section we provide a more
down to earth answer.  Exotic Lefschetz fibrations have been
constructed by fibre summing known holomorphic fibrations by
diffeomorphisms which twist the monodromy;  examples appear in the
author's thesis \cite{ivanthesis}.  These fibrations
are on manifolds which are not complex surfaces.  Examples of
non-holomorphic fibrations on manifolds homeomorphic to complex
surfaces were given by Fintushel and Stern, who distinguished the
total spaces from K\"ahler surfaces using computations of
Seiberg--Witten invariants \cite{Finsternparsin}.  Their examples were
again fibre sums.  We recall a theorem of Stipsicz \cite{Stipsicz}; a
simpler proof is given by the author in  \cite{hyperbolicdisc}:

\begin{Thm}
If a Lefschetz fibration admits a section of square $(-1)$ then it
cannot decompose as any non-trivial fibre sum.
\end{Thm}

\noindent To build exotic Lefschetz pencils, we will use a variant of
a fibre sum construction; instead of inserting a mapping class group
word of the form $\prod \delta_i = 1$ into a monodromy word, we shall
insert a more balanced word $\prod \delta_i = \prod \tau_j$ for
positive twists $\delta_i, \tau_j$.  Thus the example has the
satisfying side-effect of marrying the combinatorial and holomorphic
descriptions of Lefschetz pencils.

\begin{figure}[ht]
\begin{center}
\begin{picture}(0,0)%
\epsfig{file=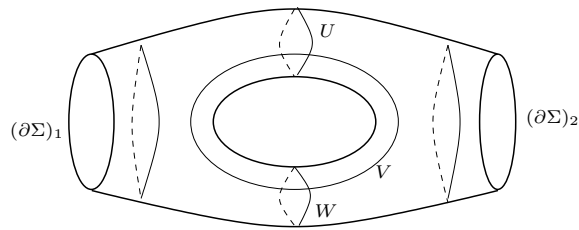}%
\end{picture}%
\setlength{\unitlength}{2368sp}%
\begingroup\makeatletter\ifx\SetFigFont\undefined%
\gdef\SetFigFont#1#2#3#4#5{%
  \reset@font\fontsize{#1}{#2pt}%
  \fontfamily{#3}\fontseries{#4}\fontshape{#5}%
  \selectfont}%
\fi\endgroup%
\begin{picture}(5400,2223)(3016,-2254)
\put(8416,-1141){\makebox(0,0)[lb]{\smash{\SetFigFont{7}{8.4}{\rmdefault}{\mddefault}{\updefault}
$(\partial \Sigma)_2$
}}}
\put(3016,-1276){\makebox(0,0)[lb]{\smash{\SetFigFont{7}{8.4}{\rmdefault}{\mddefault}{\updefault}
$(\partial \Sigma)_1$
}}}
\put(6256,-241){\makebox(0,0)[lb]{\smash{\SetFigFont{7}{8.4}{\rmdefault}{\mddefault}{\updefault}
$U$
}}}
\put(6841,-1726){\makebox(0,0)[lb]{\smash{\SetFigFont{7}{8.4}{\rmdefault}{\mddefault}{\updefault}
$V$
}}}
\put(6211,-2131){\makebox(0,0)[lb]{\smash{\SetFigFont{7}{8.4}{\rmdefault}{\mddefault}{\updefault}
$W$
}}}
\end{picture}
\end{center}
\caption{Supporting curves for a basic relation in $\Gamma_{1,2}$
  \label{relationtwists}}
\end{figure}

\begin{Example} \label{wordrelationtwists}
The mapping class group of a genus $g$ surface can be generated by
positive Dehn twists subject to relations supported in twice-holed
tori and four-punctured spheres \cite{Luo}.  In this presentation,
only one ``basic relation'' equates two non-trivial
products of positive twists.  Let $\Sigma_{1,2}$ be a torus with two
boundary circles.  Write $\delta_1, \delta_2$ for the positive twists
about curves parallel to the two boundary components $(\partial
\Sigma)_i$; there is a relation
$$(\tau_u \tau_v \tau_w)^4 \ = \ \delta_1 \delta_2,$$
where $\tau_u$ denotes the twist about the curve labelled $U$ in
Figure \ref{relationtwists}, etc.
\end{Example}

\noindent Write $\mathcal{T}$ for the operation on positive relations
which replaces a string $(\delta_1 \delta_2)$ by the string $(\tau_u
\tau_v \tau_w)^4$.  An easy computation using Novikov additivity shows
that $\mathcal{T}$ has the following effect on the topological
invariants of the four-manifold:
$$\sigma \mapsto \sigma-6;  \ e \mapsto e+10;  \ c_1 ^2 \mapsto c_1 ^2
+ 2.$$
\noindent These formulae must be modified if either of the boundary
curves in the inserted copy of $\Sigma_{1,2}$ bounds in the higher
genus surface;  for instance, the signature is changed by eight and
the Euler characteristic by twelve when fibre summing elliptic
fibrations.  Assuming we are in the generic situation, however, and
for pencils with no
reducible elements, that is discounting the $\Delta_{i \geq 1}$, we
deduce
\begin{Eqn} \label{changehits}
[\sss^2] \cdot (a \lambda - b \Delta) \
\stackrel{\mathcal{T}}{\longrightarrow} \ [\sss^2]\cdot(a\lambda -
b\Delta) - (10b-a).
\end{Eqn}
\noindent It follows that applying the operation $\mathcal{T}$ to a
Lefschetz fibration provides a clean way of decreasing the
intersection numbers with divisors.  For if $a > 11b$ then in fact the
divisor is ample and -- in the vein of the remarks
(\ref{generaltype}) -- we are unlikely to find any topological
consequence of an inclusion $\sss^2 \subset D$.  Thus we are most
likely to deal with divisors for which $a \leq 11b$; when also $a <
10b$ then the mapping class group insertion $\mathcal{T}$ will not
increase the algebraic intersection with $D$.  Moreover we have the
following key lemma:

\begin{Lem} \label{keepzerosquare}
Let $\langle \prod \gamma_i = 1 \rangle$ be a positive relation
describing a Lefschetz fibration which contains an exceptional
section.  Assume the generic fibre genus is at least two.  Then a
fibration obtained by applying either $\mathcal{T}$
or $\mathcal{T}^{-1}$ to the relation also admits an exceptional
section.
\end{Lem}

\begin{Pf}
It is enough to see that the relation $\langle (\tau_u \tau_v \tau_w)^4
\delta_2^{-1} \delta_1 ^{-1} = 1 \rangle$  describes a (non-symplectic)
fibration with a section of square zero.  For then we can view the
operations $\mathcal{T}^{\pm 1}$ as given by fibre summing two
fibrations, one with a section of square $-1$ and one with a section
of square $0$, and then excising a piece of the fibration with trivial
monodromy.  This trivial piece is determined by a relation $\langle (\prod
f_i) (\prod f_i) ^{-1} = 1 \rangle$ and such relations always admit
square zero sections.  One can then perform the fibre sums and
excisions relative to a base-point (noting for instance that every
diffeomorphism of a surface is isotopic to a once-pointed
diffeomorphism) to obtain the result.

\noindent The existence of the section of square zero for the basic
relation (\ref{wordrelationtwists}) follows from the methods of
\cite{hyperbolicdisc}.  Since the genus of the generic fibre $g \geq 2$ we
know that the supporting curves for the $\mathcal{T}$--relation do not fill the
surface.  Hence we can lift all of the individual twists to the
hyperbolic disc, the universal cover of a single smooth fibre, so as
to preserve some union of geodesics.  Look at the monodromy at the
circle at infinity which is the boundary of the disc;  after we have
lifted each (positive or negative) Dehn twist, we have lifted the
identity and have a hyperbolic automorphism of the disc.  But this
automorphism can be taken to fix a geodesic, so must
be the identity.  Since the individual twists fixed points on $\sss^1
_{\infty}$ the total rotation number of the circle at infinity under
the sequence of lifts is zero. But according to the results of
\cite{hyperbolicdisc} this precisely constructs a section of square zero.
\end{Pf}

\noindent We now return to divisors in moduli space.  The first of the
Brill--Noether divisors already introduced is the hyperelliptic divisor
$\mathcal{H}_3$ inside $\overline{M}_3$.  The cohomology class  of
this divisor in terms of the standard generators
$\lambda, \Delta_i$ for $H^2( \mgbar)$, up to a positive rational
multiple, is given by
\begin{Eqn} \label{hyperellipticlocus}
[\overline{\mathcal{H}}_3] \ = \ 9 \lambda - \Delta_0 - 3 \Delta_1.
\end{Eqn}
\noindent Technically in the above form we have given the relationship
in the Picard group of the moduli functor and not the cohomology
(Chow) ring
of the moduli space:  the discrepancy arises because of the presence
of an orbifold structure on the entire component of the stable divisor
comprising curves with elliptic tails and on the locus of
hyperelliptic curves itself. Thus translating instead to the Chow
group (and abusing notation by using the same symbols to denote the
respective generators) we have
$$9 \lambda - \Delta_0 - 3 \Delta_1 \ \ \mapsto \ \ 18 \lambda -
2\Delta_0 - 3\Delta_1;$$
\noindent in our examples $\Delta_1$ will vanish, and the re-scaling
will be inconsequential.  The particular advantage of working with the
hyperelliptic divisor is the well-known restriction on Lefschetz
pencils giving rise to spheres with image contained inside it \cite{Endo}:

\begin{Lem}[Endo]
Let $X \rightarrow \sss^2$ be the total space of a genus three
hyperelliptic Lefschetz fibration.  Then 
$$\sigma(X) \ = \ -\frac{4}{7} i + \frac{1}{7} r$$
\noindent where $i,r$ denote the numbers of irreducible and reducible
singular fibres respectively\footnote{Denoting the total number of
  singular fibres by $s$ we arrive at the taxing formula $i+r=s$.}.
\end{Lem}

\noindent There are various proofs of this result.  For topologists,
  hyperelliptic 
  fibrations are globally double branched covers of sphere bundles
  over spheres and this yields an explicit signature formula.  For
  geometers, the (nearly ample) Hodge class $\lambda$ restricts to the
  moduli space $\mathcal{H}_3$ as a certain union of boundary divisors
  since the open locus of hyperelliptic curves is affine, and the
  signature formula results from comparing to
  (\ref{signatureformula}).  Endo's original proof is an
  algebraic formulation of the last statement, analysing the signature
  cocycle in the second cohomology of the symplectic group under
  restriction to the hyperelliptic mapping class group.  In any case,
  we have a new obstruction to holomorphic structures on genus three
  fibrations: 

\begin{Prop} \label{genusthreeobstruction}
Let $X \rightarrow \sss^2$ be a genus three Lefschetz fibration with
irreducible fibres.
Suppose that 
\begin{enumerate}
\item $e(X) + 1$ is not divisible by $7$;
\item $9 \sigma(X) + 5e(X) + 40 < 0$.
\end{enumerate}

\noindent Then the Lefschetz fibration is not isotopic to a
holomorphic fibration.
\end{Prop}

\begin{Pf}
If $e+1$ is not divisible by seven then the integrality of the
signature and Endo's formula show that the fibration is not
hyperelliptic.  On the other hand, the condition that $9 \sigma(X) +
5e(X) + 40 < 0$ is exactly equivalent to the statement that
$[\sss^2_X] \cdot \mathcal{H}_3 < 0$ and so the fibration cannot yield
a sphere isotopic to a rational curve.
\end{Pf}

\noindent We now hit a catch.  There are rather few known genus three
fibrations which admit exceptional sections, and the numerology
conspires against us:  for none of these does applying the relation
$\mathcal{T}$ yield a negative algebraic intersection with the
hyperelliptic divisor.  Indeed for most we cannot apply $\mathcal{T}$
at all, since no sum of two vanishing cycles is homologically trivial,
as must be the case for the curves $U,W$ in Figure \ref{relationtwists}.
On the other hand, if we apply $\mathcal{T}^{-1}$ then we increase the
intersection with the hyperelliptic divisor.  Since the only
explicitly known relations correspond to holomorphic fibrations, this
would appear to be a losing strategy.  Fortunately, there is a
loophole.  We can start with a hyperelliptic holomorphic fibration,
for which the rational curve lies inside $\mathcal{H}_3$ and has
sufficiently negative intersection with the divisor that applying
$\mathcal{T}^{-1}$ does not destroy that property; but by an easy
count, it will necessarily destroy the hyperellipticity.  It is fairly
easy to classify holomorphic hyperelliptic fibrations with no
reducible fibres; to satisfy the constraints of
(\ref{genusthreeobstruction})  we are reduced to essentially \emph{a single
  possible example}!  In fact the $\mathcal{T}^{-1}$--substitution, on
precisely this 
positive relation, had already been derived by Terry Fuller. His
motivation was rather different; he wanted an unusual mapping class
group word to which he could apply Kirby calculus to study covering
spaces.  He kindly donates
the following manipulations.

\begin{Example}[Terry Fuller] \label{startedhere}
Let the curves $A_i, B_i, D_2, E_2$ on a genus three curve be as drawn
in Figure 
\ref{Fullergame} and write $a_i$ etc. for the positive Dehn twist
about $A_i$.  Fuller obtains the following positive relation:
\begin{Eqn} \label{Fullerrelation}
\big ( d_2 e_2 b_2 a_2 b_1 a_1 a_3 b_2 a_2 b_1 b_3 a_3 b_2 a_2 (a_1
b_1 a_2 b_2 a_3 b_3)^{10} \big ) = 1.
\end{Eqn}
\begin{figure}[ht]
\begin{center}
\begin{picture}(0,0)%
\epsfig{file=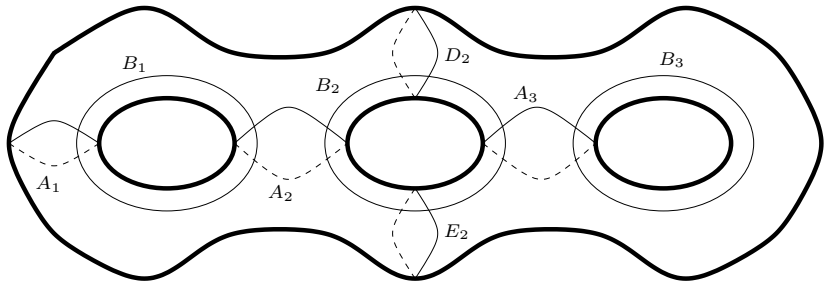}%
\end{picture}%
\setlength{\unitlength}{2368sp}%
\begingroup\makeatletter\ifx\SetFigFont\undefined%
\gdef\SetFigFont#1#2#3#4#5{%
  \reset@font\fontsize{#1}{#2pt}%
  \fontfamily{#3}\fontseries{#4}\fontshape{#5}%
  \selectfont}%
\fi\endgroup%
\begin{picture}(8188,2788)(1757,-3705)
\put(4996,-1681){\makebox(0,0)[lb]{\smash{\SetFigFont{7}{8.4}{\rmdefault}{\mddefault}{\updefault}
$B_2$
}}}
\put(2071,-2716){\makebox(0,0)[lb]{\smash{\SetFigFont{7}{8.4}{\rmdefault}{\mddefault}{\updefault}
$A_1$
}}}
\put(2971,-1456){\makebox(0,0)[lb]{\smash{\SetFigFont{7}{8.4}{\rmdefault}{\mddefault}{\updefault}
$B_1$
}}}
\put(7066,-1771){\makebox(0,0)[lb]{\smash{\SetFigFont{7}{8.4}{\rmdefault}{\mddefault}{\updefault}
$A_3$
}}}
\put(6346,-3211){\makebox(0,0)[lb]{\smash{\SetFigFont{7}{8.4}{\rmdefault}{\mddefault}{\updefault}
$E_2$
}}}
\put(4501,-2806){\makebox(0,0)[lb]{\smash{\SetFigFont{7}{8.4}{\rmdefault}{\mddefault}{\updefault}
$A_2$
}}}
\put(8596,-1411){\makebox(0,0)[lb]{\smash{\SetFigFont{7}{8.4}{\rmdefault}{\mddefault}{\updefault}
$B_3$
}}}
\put(6346,-1366){\makebox(0,0)[lb]{\smash{\SetFigFont{7}{8.4}{\rmdefault}{\mddefault}{\updefault}
$D_2$
}}}
\end{picture}
\end{center}
\caption{Supporting curves for a genus three positive relation
  \label{Fullergame}}
\end{figure}
\noindent  We begin with well-known
braid relations (which do not signify any change to the four-manifold
but only to its presentation as a positive relation):
$$a_i b_i a_i = b_i a_i b_i; \qquad a_{i+1} b_i a_{i+1} = b_i a_{i+1}
b_i.$$
\noindent Using these freely, along with the fact that Dehn
twists about disjoint curves commute, Fuller shows:
\begin{Eqn} \label{stageone}
(a_1 b_1 a_2 b_2 a_3 b_3)^2 \ = \ (a_1 b_1 a_2)^2 b_2 a_2 a_3 b_2 b_3
a_3.
\end{Eqn}
\noindent Using (\ref{stageone}), and writing 
$$(a_1 b_1 a_2 b_2 a_3 b_3)^4 \ = \ ((a_1 b_1 a_2 b_2 a_3 b_3)^2)^2$$
\noindent a further sequence of braid manipulations brings you to
\begin{Eqn} \label{stagetwo}
(a_1 b_1 a_2 b_2 a_3 b_3)^4 \ = \ (a_1 b_1 a_2)^4 b_2 a_2 b_1 a_1 a_3
b_2 a_2 b_1 b_3 a_3 b_2 a_2.
\end{Eqn}
\noindent Now employ the $\mathcal{T}^{-1}$ substitution to the
relation given below:
$$(a_1 b_1 a_2)^4 \ = \ d_2 e_2; \qquad (a_1 b_1 a_2 b_2 a_3 b_3)^{14}
= 1$$
\noindent and combine together (\ref{stageone}) and (\ref{stagetwo}):
\begin{displaymath}
\begin{array}{rl}

(a_1 b_1 a_2 b_2 a_3 b_3)^{14} = & (a_1 b_1 a_2 b_2 a_3 b_3)^4 (a_1
b_1 a_2 b_2 a_3 b_3)^{10} \\
= & (a_1 b_1 a_2)^4 b_2 a_2 b_1 a_1 a_3 b_2 a_2 b_1 b_3 a_3 b_2 a_2 (a_1
b_1 a_2 b_2 a_3 b_3)^{10} \\
 = & d_2 e_2 b_2 a_2 b_1 a_1 a_3 b_2 a_2 b_1 b_3 a_3 b_2 a_2 (a_1
b_1 a_2 b_2 a_3 b_3)^{10}.

\end{array}
\end{displaymath}
\noindent This gives us the relation we require.
\end{Example}

\noindent The relation $\langle (a_1 \ldots b_3)^{14} = 1 \rangle$
defines the monodromy of a hyperelliptic surface obtained by taking
the double branched cover of $\ff_2$ over a curve in the class $7|s_0|
\amalg |s_{\infty}|$, where $s_0, s_{\infty}$ denote the sections of
square $\pm 2$ respectively.  Since the section of square two lies
inside the branch locus, it lifts to an exceptional section of the
genus three fibration;  if we blow this down we obtain a simply
connected complex surface of general type.  To prove that this does indeed
give the genus three fibration we require, consider the curve in $|7
s_0|$ given locally by $z_1 ^7 + z_2 ^{14} = 0$.  Deform the
singularity to a union of $14$ simple tangencies between distinct
sheets of the curve.  At each of these points the new equation is
equivalent to $u^2 + v^7 = 0$ in suitable holomorphic co-ordinates,
which gives the monodromy of the $(2,7)$--torus knot (which is fibred
of genus three).  Then compute that the monodromy of this knot is
indeed given by the product $(a_1 b_1 a_2 b_2 a_3 b_3)$.

\noindent Write $W$ for the total space of the Lefschetz fibration
defined by (\ref{startedhere}).  The pencil of curves, both
before and after the modification by $\mathcal{T}^{-1}$, has only
irreducible fibres and hence does define a canonical symplectic
structure by the theorem of Gompf mentioned in the Introduction.  The
number of singular fibres in the modified pencil is $74$ and hence
$e(W) = 66$.  Moreover the signature of the Horikawa surface is $-48$
(by Endo's formula, say) and so the modified manifold has $\sigma(W) =
-42$. (This checks with the result obtained, via a computer-implemented
algorithm based on Wall's non-additivity, by Ozbagci in
\cite{Ozbagci}.  In particular the signature is computable from the
mapping class group word, even if the derivation of that word is not
available!)  It is now a triviality to 
apply (\ref{genusthreeobstruction}) 
and deduce that we have indeed obtained a symplectic non-K\"ahler
pencil.  This completes the proof of (\ref{foundone}).

\begin{Rmk}
$W$ is homeomorphic to a simply connected complex
surface.  To the author's knowledge, the only way to prove that it might
not be K\"ahler is to compute instanton or Seiberg--Witten invariants.
Such computations seem intractable for manifolds presented via
positive relations;  note that we have deduced that the Lefschetz
structure is not holomorphic without needing to determine whether $W$
is diffeomorphic to a K\"ahler surface.
\end{Rmk}

\noindent We remark for completeness that we also find large classes
of non-holomorphic genus three Lefschetz fibrations on manifolds
homeomorphic to complex surfaces.  These examples are elementary in
the sense that they do not rely on theorems from gauge theory or
Donaldson's construction.

\begin{Example}
A hyperelliptic
genus three Lefschetz fibration $Z$ with $\delta$ singular fibres, all
irreducible, has $\delta = 7r$ for some integer $r$ and topological
invariants $e = 7r-8, \, \sigma = -4r$.  Now under a fibre summation the
signature and Euler characteristic transform as
$$\sigma \ = \ \sigma_1 + \sigma_2; \qquad e \ = \ e_1 + e_2 - 2e(F).$$
\noindent It follows immediately that $Z \sharp_{\Sigma_3}
W$ will satisfy (\ref{genusthreeobstruction}) and also the
Noether inequality:
$$(9 \sigma + 5e + 40)(Z \sharp_F W) \ = \ (9 \sigma + 5e + 40)(W) -
r$$
$$(5c_1 ^2 - c_2 +36)(Z \sharp_F W) \ = \ (5c_1^2 -c_2 +36)(W) + 3r.$$
\noindent One can also check that for any $Z$ there is a minimal
complex surface of general type homeomorphic to $Z \sharp_F W$ using the
geographic criterion of (\cite{BPV}, VII (9.1)).  Indeed for these
manifolds it is known that there is a simply connected such
complex surface \cite{PerssonComp}.
By similar
manipulations the reader can check the relevant
numbers in the case of fibre summing $W$ with itself.  Note that many
simply-connected hyperelliptic genus three fibrations are K\"ahler,
but the holomorphic 
structure of the fibrations is lost on summing with $W$.
\end{Example}

\noindent By a result of Stipsicz \cite{Stipsicz2} we know
that many of these manifolds are minimal.  Auroux (private
conversation) has noted that one can improve
the first condition in (\ref{genusthreeobstruction}) to demanding only
that $e(X) + 1$ be not
divisible by $14$, by studying the braid factorisation for a
hyperelliptic fibration instead of just the mapping class group
factorisation.

\begin{Rmk}
Suppose we have a Lefschetz pencil which is not holomorphic, for
instance by arguments as above.  If the canonical symplectic form is
isotopic to a K\"ahler form, then we can use the Riemann--Roch theorem
to estimate the number of sections of the line bundle $\mathcal{L}_C$
with first Chern class dual to the fibre $C$ of the pencil.  If any
section of this line bundle had smooth zero set, then a generic pencil
of sections would define a Lefschetz pencil by Bertini's theorem and
general position arguments.  The existence of such a pencil may well
be ruled out by the topological constraints on the homology class of
the sphere, which is determined by the global topology of the
manifold.  However, in general the line bundle $\mathcal{L}_C$ is
ample and not necessarily very ample, and it seems impossible (even in
examples) to prove directly that if there is a K\"ahler structure on
the manifold then there would be some smooth complex curve in the
linear system $|C|$.  This prevents us from answering our original
Question (2) in (\ref{Question}).

\noindent The situation for Question (3) is even worse:  we 
expect the homology classes realised by our pencils to
contain no smooth complex curve at all. In any complex surface with a
homology class with 
an immersed but not embedded holomorphic representative, we can smooth
nodes to build symplectic submanifolds which are uninterestingly
distinct from complex curves.
\end{Rmk}

\noindent The results of this section can be generalised to other
divisors, for instance the divisor of trigonal curves in
$\overline{M}_5$, but the author is not aware of further applications
or phenomena.


\section{Weierstrass divisors and sections}

We give a final application of (\ref{principle}) with a result
on \emph{sections} of fibrations that one can prove in a similar
vein.   This adds to the now substantial body of
knowledge on genus two fibrations \cite{ivanhodge}.
Recall the symplectic isotopy conjecture due to Siebert and Tian
\cite{ST} states that every connected
symplectic submanifold of a relatively minimal sphere bundle over a
sphere is isotopic to a complex submanifold\footnote{In lectures in
  June 2001, Siebert and Tian have announced a partial resolution of
  this conjecture which is probably strong enough to show that genus
  two Lefschetz fibrations without reducible fibres are holomorphic.
  This would suffice for this Proposition.}. 

\begin{Prop} \label{squareattwo}
Let $f\co X \rightarrow \sss^2$ be a genus two Lefschetz fibration and $s\co 
\sss^2 \rightarrow X$ a distinguished section of $f$.  Suppose the
fibration has $\delta_0$ irreducible singular fibres, $\delta_1$
reducible ones and $\delta_0 + 2\delta_1 = 10m$ for $m \in \zz$.
Suppose also that the symplectic isotopy conjecture is valid.  Then
\begin{enumerate}
\item $3 | s\cdot s | \geq m+ \delta_1$ or
\item $4 | s \cdot s| = m + \delta_1$.
\end{enumerate}
\end{Prop}

\begin{Pf} 
Inside $\overline{M}_2^1$, the moduli space of stable
curves with a single marked point -- equivalently the universal curve
$\overline{\mathcal{C}}_2$ over $\overline{M}_2$ -- there is a divisor
$\mathcal{W}$ comprising the closure of pairs $(C,p)$ where $p$ is one
of the Weierstrass points of $C$.  In the notation of
\cite{moduliofcurves}, the cohomology of
$\overline{M}_2^1$ is generated by elements $\omega_{RD}, \lambda,
\delta_0, \delta_1$, where $\omega_{RD}$ is the relative dualising sheaf
and $\lambda$ the pullback of the Hodge bundle under
$\overline{M}_2^1 \rightarrow \overline{M}_2$.  The adjunction formula
allows us to
identify the value of the class $\omega_{RD}$ on our family with the
negative of the self-intersection of the distinguished section:
\begin{Eqn} \label{RDisself-int}
\langle \omega_{RD} , [B] \rangle \ = \ - s \cdot s.
\end{Eqn}
\noindent With respect to this basis, we have an identity in the
Picard group of the moduli functor
\begin{Eqn} \label{weierstrass}
[ \mathcal{W} ] \ = \ 3\omega_{RD} - \lambda -\delta_1.
\end{Eqn}
\noindent Let $f\co  X \rightarrow \pp^1$ be a
genus two Lefschetz fibration with a distinguished section $s\co  \pp^1
\rightarrow X$.  Recall from \cite{ivanhodge} that we can write $X$ as
a double cover of a rational ruled manifold $\mathrm{Rat}_X$ over a
locus which is
smooth away from finitely many infinitely close triple points (their
number given by the number of reducible fibres in the Lefschetz
fibration).  Moreover, on any fibre the branch covering map can be
identified with the hyperelliptic involution on the fibre after
choosing metrics:  that is, the ramification locus is precisely the
closure of the union of the Weierstrass points in each fibre.

\noindent Via the fixed section $s$ we induce a map $\sss^2 \rightarrow
\overline{M}_2^1$.  If the map has image inside the divisor $\mathcal{W}$ then
the section lies inside the locus of Weierstrass points and gives rise
to a branch locus which is disconnected.   Siebert
and Tian \cite{ST} analyse the
branch locus for any hyperelliptic fibration branched over a rational
ruled manifold and
show that it has at most two components, and if it is not connected
then one of the components is a sphere section for the natural
fibration of the ruled surface over $\sss^2$.  In this case
the self-intersection of the sphere is given by $-k$ where the base
can be identified with $\pp( \mathcal{O} \oplus \mathcal{O}(k))$.
(Warning: Note here that $k$ is even and $-k$ is the self-intersection
of the sphere as a submanifold of the base;  the natural lift of the
sphere to the ramification locus has square $-k/2$.)

\noindent  On any Lefschetz fibration, there is an operation which
removes a reducible fibre and replaces it by a sequence of
$(4h+2)2h$ irreducible fibres, where $h$ is the genus of one component
of the reducible fibre. Thus for a genus two fibration we can remove a
reducible fibre and replace it by twelve irreducible ones.  This
operation can be localised
downstairs in the branched covers and corresponds to resolving or
deforming an infinitely close triple point singularity.  It follows
from this description that if there is a section disjoint from the
singularity, one can trade the two local models without changing the
section or its square.  Thus we have an operation on genus two
fibrations which has the numeric effect
$$\delta_1 \mapsto \delta_1 -1 ; \ \delta_0 \mapsto \delta_0 +12 ; \ m
\mapsto m+1; \ | s \cdot s| \mapsto |s \cdot s|.$$
\noindent If we trade all the reducible fibres for irreducible ones,
we arrive at a fibration which is (modulo the symplectic isotopy
conjecture for surfaces of appropriate bidegree) necessarily K\"ahler.

\noindent In this case, the
rational curve it defines in $\overline{M}_2 ^1$ is either contained
in $\mathcal{W}$ or has locally positive intersections with it.  Since $3 |s
\cdot s| - m - \delta_1$ is unchanged by the removal of the reducible
fibres, we  obtain the two cases of the proposition.  For suppose this
value is negative.  Then the rational curve lies inside $\mathcal{W}$
and defines a hyperelliptic fibration with disconnected branch locus.
Then by Siebert--Tian the self-intersection of the section component of
the fibration is precisely $-k/2$ where the base of the hyperelliptic
double cover is $\pp( \mathcal{O} \oplus \mathcal{O}(k))$.  Moreover  by a
result of \cite{ivanhodge}, we can relate this value $k$ to the number
of critical fibres of the fibration: precisely, $|k| = m$ where there
are $10m$ critical fibres.

\noindent Since we assume the symplectic isotopy conjecture, a genus
two Lefschetz fibration with only irreducible fibres is holomorphic,
and such objects were classified by Chakiris and 
independently by the author \cite{ivanhodge}.  The only irreducible
fibrations are listed in the following proposition
(\ref{complexlist}).  Using the
classification, and the explicit constructions of \cite{ivanhodge},
one can check that if the ramification locus contains
a section component, then its square is related to the number of
critical fibres $10m$ by the formula $|s \cdot s| = m/4$.
Since after trading away the reducible
fibres this number $m$ is
given by $m+\delta_1$ for the original fibration, we find in the end
that $4 |s \cdot s| = m+\delta_1$ as claimed.
\end{Pf}

\noindent One consequence of this result is a strong restriction on
which genus two Lefschetz fibrations can admit sections of square
$(-1)$.  This corollary can be obtained by other methods,
which may be of interest in themselves.  In Figure \ref{Fullergame} we
depict a genus three curve: if we cut this along the cycles $D_2$ and
$E_2$ and glue the boundary curves of the left component, we obtain a
genus two curve with distinguished cycles $A_1, B_1, A_2, B_2$ and
$D_2=E_2=A_3$. Again we will use lower case letters to denote the associated
Dehn twist diffeomorphisms, now of the genus two curve.  The
idea of this second proof is that it is easier to
restrict the topology of Lefschetz pencils than Lefschetz fibrations.
Recall that we showed in (\ref{notalwayseasy}) that the Poincar\'e dual of the
fibre $[C] \in H_2 (X; \zz)$ of a pencil was not always represented by
a symplectic form 
adapted to the fibration.  If there were such a form, and if the
four-manifold satisfied the constraint $b_+ > 1$, then by Taubes'
results on the canonical class we would deduce that $K_X \cdot [C]
\geq 0$ with equality if and only if $K_X = 0$.  This helpful
property persists:

\begin{Prop} \label{likesymplectic}
Let $X$ be a smooth four-manifold with a Lefschetz pencil of curves
each representing a homology class $[C]$.  Suppose that $b_+ (X) >
1$.  Then the canonical
deformation equivalence class of symplectic forms on $X$ defined by
the pencil gives an almost complex structure such that $K_X$ satisfies
$$K_X \cdot [C] \geq 0; \ \qquad \ K_X \cdot [C] = 0 \ \Rightarrow \
K_X = 0.$$
\end{Prop}

\begin{Pf}
The canonical class for the associated Lefschetz fibration represents
$\pi^* K_X + \sum E_i$ in homology, where the $E_i$ are the
exceptional curves and $\pi\co  X' \rightarrow X$ the blow-down map.  By
adjunction, we know that $K_{X'} \cdot [\pi^* C] = \deg (K_C) =
2g-2$.  Moreover we know that for any almost complex structure $J$,
there is an immersed holomorphic representative for $K_X$ and that
pseudoholomorphic curves have locally positive intersections.
Choosing $J$ so that the exceptional spheres are indeed
pseudoholomorphic, it
follows that $[C] \cdot [C] 
\leq 2g-2$ with equality iff $K_{X'} = \sum E_i$.  But this implies
the proposition.
\end{Pf}

\noindent Armed with this we can establish some control on the homotopy types
of manifolds with genus two pencils\footnote{A stronger version of
  this result follows from combining the isotopy conjecture and the
  first result of the section.}.

\begin{Thm} \label{complexlist}
Let a symplectic four-manifold $X$ admit a
Lefschetz pencil of genus two curves.  Then there are only finitely
many possibilities for the numbers $n,s$ of irreducible and reducible
critical fibres, or equivalently for the pair $(c_1 ^2 (X), c_2(X))$.
In particular $e(X) < 40$.  If $s=0$ then the fibration determined by
the pencil on $X$ is
homeomorphic to one of the
simply connected complex surfaces  associated to the three monodromy words:

\begin{enumerate}
\item $(a_1 b_1 a_2 b_2 a_3 a_3 b_2 a_2 b_1 a_1)^2 = 1$;
\item $(a_1 b_1 a_2 b_2 a_3)^6 = 1$;
\item $(a_1 b_1 a_2 b_2)^{10} = 1$.
\end{enumerate}
\end{Thm}

\begin{Pf}
Be given $X$ and write $[C]$ for the homology class defined by a fibre
  of the pencil.
  Assume  $b_+ (X) > 1$; then our proposition (\ref{likesymplectic})
  tells us that $2g-2 = 2 = K_X \cdot [C] + [C]^2$
with both of the terms in the last expression non-negative, and the
  first zero only if $K_X = 0$.  This
 gives a small number of possibilities:  either $K_X = 0$ and
 $\omega^2 = 2$, or $K_X
 \cdot \omega = 1$ and $\omega^2 = 1$.  Suppose first that $K_X =
 0$;  then $c_1 ^2 = 2e+3\sigma = 0$.  The Euler characteristic
 $e$ and signature $\sigma$ are
 determined by the numbers of critical fibres.  Let there be $n$
 non-separating critical fibres and $s$ separating ones.  Then
$$e = n+s-6; \qquad \sigma = \frac{3}{5} n - \frac{1}{5}s +2$$
\noindent where we have used the formulae for the associated fibration
of curves from \cite{ivanhodge} and the fact that the pencil has two
base-points.  It follows from these formulae that 
$$n+7s = 30$$
\noindent whilst (for any genus two fibration, since
$(\Gamma_2)_{\mathrm{ab}} = \zz_{10}$) also $n+2s$ is divisible by
$10$.  This gives two 
possibilities: $n=30, s=0$ and $n=16,s=2$.  In the first case, we know
the fibration is simply-connected \cite{ST} and (since $K_X = 0$)
it is minimal.  An easy computation then
shows that it is homeomorphic to the $K3$ surface blown up twice,
which is described by the second word listed in the statement of the
proposition.  In fact a symplectic
four-manifold with $\pi_1 = 0, c_1 = 0$ is necessarily homeomorphic
to the $K3$ surface by a result of Morgan and Szabo
\cite{Morg-Sz}.

\noindent Suppose instead that
$\omega^2 = 1$.  Let  $X' = X
\sharp \overline{\cc \pp}^2$ be the
total space of the Lefschetz fibration.  By Taubes, we have an
immersed holomorphic curve representing $K_{X'}$ and containing the
exceptional section.  
If the
symplectic 
representative for $K_{X'}$ is smooth, then it meets each genus two
fibre two times; since it contains the exceptional section, it must
comprise two disjoint sections.  But by adjunction this gives a 
smooth (symplectic) genus zero representative for $K_X$, which is
therefore an exceptional sphere.  Thus $X$ is the blow-up of a simply
connected manifold with $K_X = 0$, and by the result of
Morgan--Szabo referred to above we see that the given pencil of
curves is just the blow-up of the usual genus two pencil
on $K3$ at one of its two base-points.  In particular, the associated
fibration is the one obtained above.

\noindent The only other possibility is that the symplectic subvariety
provided by Taubes is in fact composed of several components.  This
subvariety meets every genus two fibre with which it shares no
component locally positively and with algebraic intersection number
two.  Hence the curve must comprise the
exceptional section counted to multiplicity two, and a number $r$ of
fibres of the fibration.  Suppose we write $K_{X'} = 2E + rF$;  then $K_{X'}^2
= 4r-4$ and so $K_X ^2 = 4r-3$, blowing down again.  However, this is
also given by $g-1$ where $g$ is the genus of the smooth symplectic
representative for $K_X$ obtained by Taubes for a generic complex
structure.  (Here we define genus via a sum over components if
necessary).  By construction, one symplectic representative for
$K_X$ is given by a number $r$ of the genus two curves of the pencil
smoothed at the base-point.  The result is that
$$2r = 4r-2$$
\noindent and hence $r=1$.  But then this determines $c_1 ^2 (X') = 0$
and this is enough to fix the number of critical fibres of the
fibration, using the usual formulae for $e, \sigma$.  This gives the
finiteness we require.   If we also know $s=0$ and 
the fibration is simply connected, and since the intersection form
must be odd as $(K_X) ^2 = 1$, we have determined $X$ to
homeomorphism.  It must be a surface of general type, with $K_X =
[\omega]$ represented by a genus two curve of square one.  Such is
described by the third word in the list of monodromies given in the
statement of the theorem.

\noindent  The last remaining case is where Taubes does not apply,
that is $b_+ (X) = 1$.  In this case 
$$\sigma \ = \ 1-b_- \ = \ \frac{3}{5}n - \frac{1}{5}s + \omega^2$$
\noindent and 
$$e \ = \ n+s - 4 - \omega^2 \ = \ 3-2b_1 + b_-.$$
\noindent Moreover we still know that $K_X \cdot \omega + \omega^2 =
2$ and $b_1 = 0$ or $b_1 = 2$, since $b_1, b_+$ have opposite parity
on any almost complex four-manifold and for any Lefschetz fibration
$b_1$ of the total space is strictly smaller than $b_1$ of the fibre.
For a fibration with $s=0$ we know we have simply connected total
space, so $b_1 = 0$ and some easy manipulations give $\sigma + e = 4$
and then $n=20$.  To homeomorphism this
gives the manifold given by the first monodromy word on the
list above.  The
other cases, with $s \neq 0$ and either $b_1$ are also easily listed:
in all cases $n \leq 20$ and hence only finitely many pairs $(n,s)$ arise.
\end{Pf}

\noindent This establishes the last of the theorems given in the
opening section of the paper.  Although determining a symplectic
manifold only to homeomorphism is a very weak statement, the limited
geography of manifolds with genus two pencils is striking
in its own right.  It would be very interesting to
understand if there is any analogue of this at higher genera.



\begin{thebibliography}

\bibitem{arabel} {\bf E~Arabello}, {\bf M~Cornalba}, \emph{The
{P}icard groups of the moduli spaces of curves}, Topology, {26} (1987)
153--171

\bibitem{Auroux} {\bf D~Auroux}, \emph{Symplectic four-manifolds as
branched coverings over {$\cc \pp^2$}}, Invent. Math. {139} (2000)
551--602

\bibitem{AKstabilisation} {\bf D~Auroux}, {\bf L~Katzarkov}, \emph{The
degree doubling formula for braid monodromies and {L}efschetz
pencils}, Preprint (2000)

\bibitem{BPV} {\bf W~Barth}, {\bf C~Peters}, {\bf A Van\,de\,Ven},
\emph{Compact complex surfaces}, Springer (1984)

\bibitem{SKDIS} {\bf S\,K~Donaldson}, {\bf I~Smith}, \emph{Lefschetz
pencils and the canonical class for symplectic $4$--manifolds},
Preprint (2000) {\tt arxiv:math.SG/0012067}

\bibitem{Donaldson:pencils} {\bf S\,K Donaldson}, \emph{Lefschetz
pencils on symplectic manifolds}, J. Diff.  Geom. {53} (1999) 205--236

\bibitem{Endo} {\bf H~Endo}, \emph{Meyer's signature cocyle and
hyperelliptic fibrations}, Math.  Ann. {316} (2000) 237--257

\bibitem{Finsternparsin} {\bf R~Fintushel}, {\bf R~Stern},
\emph{Counterexamples to a symplectic analogue of a theorem of
{Arakelov and Parsin}}, Preprint (1999)

\bibitem{GompfGokova} {\bf R~Gompf}, \emph{The topology of symplectic
manifolds}, Turkish J. Math.  {25} (2001) 43--59

\bibitem{Harer} {\bf J~Harer}, \emph{The second homology group of the
mapping class group of an orientable surface}, Invent. Math. {72}
(1983) 221--239

\bibitem{moduliofcurves} {\bf J~Harris}, {\bf I~Morrison},
\emph{Moduli of curves}, Springer (1998)

\bibitem{LiuSW} {\bf A~Liu}, \emph{Some new applications of the
general wall-crossing formula}, Math. Res. Lett. {3} (1996) 569--585

\bibitem{Luo} {\bf F~Luo}, \emph{A presentation of the mapping class
groups}, Math. Res. Lett.  {4} (1997) 735--739

\bibitem{Matsumoto} {\bf Y~Matsumoto}, \emph{Lefschetz fibrations of
genus two -- a topological approach}, The 37th Taniguchi Symposium on
topology and Teichm\"uller spaces (S~Kojima et~al, eds.) World
Scientific (1996)

\bibitem{McDS2} {\bf D~McDuff}, {\bf D~Salamon}, \emph{A survey of
symplectic $4$--manifolds with $b_+ = 1$}, Turkish J. Math. {20}
(1996) 47--61

\bibitem{Morg-Sz} {\bf J~Morgan}, {\bf Z~Szabo}, \emph{Homotopy {$K3$}
surfaces and mod $2$ {Seiberg--Witten} invariants},
Math. Res. Lett. {4} (1997) 17--21

\bibitem{Ozbagci} {\bf B~Ozbagci}, \emph{Signatures of {L}efschetz
fibrations}, Pacific J. Math. (to appear)

\bibitem{PerssonComp} {\bf U~Persson}, \emph{Chern invariants of
surfaces of general type}, Comp. Math.  {43} (1981) 3--58

\bibitem{ST} {\bf B~Siebert}, {\bf G~Tian}, \emph{On hyperelliptic
{$C^{\infty}$--Lefschetz} fibrations of four-manifolds},
Commun. Contemp. Math. {1} (1999) 466--488

\bibitem{ivaninstanton} {\bf I~Smith}, \emph{Gauge theory and
symplectic linear systems}, In preparation

\bibitem{hyperbolicdisc} {\bf I~Smith}, \emph{Geometric monodromy and
the hyperbolic disc}, Quarterly J. Math. (Oxford) 52 (2001) 217--228

\bibitem{ivanSDforSS} {\bf I~Smith}, \emph{{Serre--Taubes} duality for
pseudoholomorphic curves}, Preprint (2001)

\bibitem{ivanthesis} {\bf I~Smith}, \emph{Symplectic geometry of
{L}efschetz fibrations}, Ph.D. thesis, Oxford University (1998)

\bibitem{ivanhodge} {\bf I~Smith}, \emph{Lefschetz fibrations and the
{H}odge bundle}, Geometry and Topology {3} (1999) 211--233

\bibitem{Stipsicz2} {\bf A~Stipsicz}, \emph{On the number of vanishing
cycles in {L}efschetz fibrations}, Math. Res. Lett. {6} (1999)
449--456

\bibitem{Stipsicz} {\bf A~Stipsicz}, \emph{On the indecomposability of
certain {Lefschetz} fibrations}, Proc. Amer. Math. Soc. {129} (2001)
1499--1502

\bibitem{Taubes} {\bf C\,H Taubes}, \emph{The {S}eiberg--{W}itten and
the {G}romov invariants}, Math.  Res. Letters, {2} (1995) 221--238

\end{thebibliography}
\end{document}